\documentclass[12pt]{article}

\usepackage{geometry}                % See geometry.pdf to learn the layout options. There are lots.
\usepackage{graphicx}
\usepackage{amssymb,amsmath}
\usepackage{epstopdf}
\usepackage[leqno,fleqn,intlimits]{empheq}
\renewcommand{\i}{\boldsymbol{i}} % the imaginary quaternions
\renewcommand{\j}{\boldsymbol{j}}
\renewcommand{\k}{\boldsymbol{k}} 

\newcommand{\I}{{\scriptscriptstyle I}} % makes the script smaller
\newcommand{\J}{{\scriptscriptstyle J}}
\newcommand{\K}{{\scriptscriptstyle K}} 
\newcommand{\0}{{\scriptstyle 0}}
\newcommand{\1}{{\scriptstyle 1}}
\newcommand{\2}{{\scriptstyle 2}}

\newcommand{\sprod}{\!\cdot\!}   % scalar product
\newcommand{\vprod}{\!\times\!} % vector product

\newcommand{\V}{V} % the potential

\newcommand{\pt}{\hspace{1pt}} % really small space

\title{Quaternionic K\"ahler spaces with large toric symmetry}
\author{Radu A. Iona\c{s}\,\footnote{Email: ionas@max2.physics.sunysb.edu}} 
\date{}                                           % Activate to display a given date or no date

\begin{document}
\maketitle
\thispagestyle{empty}
\setcounter{page}{0}
\setcounter{section}{-1}

\begin{abstract}
We consider a general $4n$-dimensional quaternionic K\"ahler geometry with a free action of the torus $T^{n+1}$. The toric action lifts onto the Swann bundle of the quaternionic K\"ahler space to a tri-holomorphic action that commutes with the standard $\mathbb{H}^{\times}$ action on the bundle. By matching Pedersen and Poon's generalized Gibbons-Hawking Ansatz description of the total space with the Swann picture we extract the local geometry of the  quaternionic K\"ahler base. Specifically, we obtain explicit expressions for the quaternionic K\"ahler metric and $Sp(1)$ connection in terms of a set of reduced Higgs fields and connection 1-forms that satisfy a reduced Bogomol'nyi-type equation. We find, moreover, that these Higgs fields can be derived from a single function $\V$ satisfying a system of linear second-order partial differential constraints. In four dimensions, corresponding to the case of self-dual Einstein manifolds with two commuting Killing vector fields, our formulas coincide with those obtained through a different approach by Calderbank and Pedersen. Finally, we show how to construct explicit solutions to the reduced Bogomol'nyi and $\V$ equations by means of Lindstr\"om and Ro\v{c}ek's generalized Legendre transform construction for a large class of quaternionic K\"ahler manifolds related to the c-map.
\end{abstract}

\newpage

\section{Introduction}

Quaternionic K\"ahler spaces with large toric actions have been long studied, from various perspectives, in both physics and mathematics. On the mathematical side, one of the most interesting results in this direction was obtained in four dimensions by Calderbank and Pedersen \cite{MR1950174}, who gave a complete local classification of self-dual Einstein metrics of non-zero scalar curvature with two commuting Killing vector fields; any such metric was shown to have an explicit local form determined by an eigenvalue of the Laplacian on the hyperbolic plane. This was achieved by an intrinsic, four-dimensional approach, through an ingenious  interplay among Joyce's classification of self-dual manifolds with two commuting surface-orthogonal conformal vector fields \cite{MR1324633}, Tod's description of self-dual Einstein metrics with non-zero scalar curvature and admitting a Killing field in terms of the $SU(\infty)$ Toda equation \cite{MR1423177,MR1097787} and Ward's theory of axisymmetric Einstein-Weyl spaces \cite{MR1045295}. The corresponding Swann bundles -- hyperk\"ahler spaces of one quaternionic dimension higher whose hyperk\"ahler structure encodes the quaternionic K\"ahler structure of the base -- were constructed more as an afterthought, through a bottom-up approach.

In physics, quaternionic K\"ahler manifolds arise usually as target spaces of locally supersymmetric sigma models with eight supercharges \cite{Bagger:1983tt}. They describe the hypermultiplet moduli spaces of type II strings compactified on Calabi-Yau 3-folds \cite{Cecotti:1988qn} or heterotic strings compactified on K3 surfaces \cite{Aspinwall:1998bw,Witten:1999fq}. Swann bundles, on the other hand, arise as target spaces of field theories that are invariant under rigid $N=2$ superconformal symmetry \cite{deWit:1999fp,deWit:2001dj}. The process through which one retrieves a quaternionic K\"ahler manifold from its Swann bundle is known as $N=2$ superconformal quotient. In particular, quaternionic K\"ahler spaces of the type that we consider here occur for example in relation to the c-map \cite{Cecotti:1988qn,Ferrara:1989ik}. This is a construction which can be understood in the context of T-duality between type IIA and type IIB string theories compactified on circles of inverse radii, and which maps projective special K\"ahler manifolds of complex dimension $n-1$ to a certain class of quaternionic K\"ahler manifolds of quaternionic dimension $n$, admitting, among other symmetries, a set of $n+1$ commuting Killing vector fields. Due to its association to special K\"ahler geometry, the c-map has been discussed in connection to topological strings and,  by way of the Ooguri-Vafa-Strominger conjecture \cite{Ooguri:2004zv}, to the Bekenstein-Hawking entropy of supersymmetric black holes \cite{Rocek:2005ij,Rocek:2006xb,Neitzke:2007ke}.

A field-theoretic rederivation of the Calderbank-Pedersen metric was given in \cite{deWit:2006gn}. This analysis was later extended to the case of eight dimensions in \cite{deWit:2007qz}. The superconformal quotient of arbitrary $4n$-dimensional quaternionic K\"ahler manifolds with $n+1$ commuting Killing vector fields was also considered in \cite{Alexandrov:2008nk} from a different point of view than the one we take here.

This paper is organized as follows: In section 1 we recall a few basic facts about quaternionic K\"ahler manifolds and Swann bundles. In sections 2 and 3 we review Pedersen and Poon's generalization of the Gibbons-Hawking Ansatz and Lindstr\"om and Ro\v{c}ek's Legendre transform approach to constructing toric hyperk\"ahler varieties from twistor lines and meromorphic functions. In particular, we discuss the Legendre transform construction of toric Swann bundles, the constraints one must impose in this case on the meromorphic functions, and explain the relation between the collective degrees of freedom of  the associated monopole configurations and the $\mathbb{H}^{\times}$ action orbits. In section 4 we use the equivariance of the moment maps associated to the toric action on the Swann bundle to coordinatize the quaternionic K\"ahler base. Dissecting the generalized Gibbons-Hawking Higgs fields and connection 1-forms along the lines of the Swann fiber structure allows us then to project the abelian monopole equations onto the base and derive eventually closed-form expressions for the quaternionic K\"ahler metric and $Sp(1)$ connection. These are given in terms of a set of reduced Higgs fields and connections; the latter satisfy a Bogomol'nyi-type equation while the former turn out to be determined by a single real-valued function closely related to the hyperk\"ahler potential of the Swann bundle. In section 5 we show how these results match those of Calderbank and Pedersen when particularized to dimension four. In section 6 we give a twistor-theoretic prescription on how to construct explicit solutions to the previously obtained equations for the class of quaternionic K\"ahler manifolds given by the c-map and then work out  a few examples.

\section{Quaternionic K\"ahler manifolds and Swann \\ bundles} \label{SEC:QK-&-SB}

A $4n$-dimensional Riemannian manifold ${\cal M}$ is called quaternionic K\"ahler if its holonomy group is contained in the $Sp(n)Sp(1)$ subgroup of $SO(4n)$.\footnote{ In this paper we will be considering pseudo-quaternionic K\"ahler manifolds as well, in which case this definition has to be modified accordingly, see {\it e.g.} \cite{MR2140713}.} This is a non-trivial restriction for $n>1$, but for $n=1$ it trivially yields all oriented four-dimensional Riemannian manifolds, since $Sp(1)Sp(1) \simeq SO(4)$. The precise analogue of quaternionic K\"ahler manifolds in dimension $4$ are the Einstein self-dual manifolds. 

Manifolds with the quaternionic K\"ahler property possess a rank-3 subbundle of $\mbox{End\,}T{\cal M}$ spanned  locally by a basis formed of a triplet  of almost complex structures $I$, $J$, $K$ that satisfy the algebra of quaternions ($I^2=J^2=K^2=-1$, $IJ=K$, a.s.o.) as well as a metric $g$ which is Hermitian with respect to each of the almost complex structures. From these objects one can further construct a bundle of almost K\"ahler structures $\Theta_1$, $\Theta_2$, $\Theta_3$ from $\wedge^2T^*{\cal M}$ by taking $\Theta_1(X,Y) = g(IX,Y)$ a.s.o., for all $X,Y \in T{\cal M}$; it is convenient to assemble them into an $ \mbox{Im\,}\mathbb{H}$-valued 2-form, $\Theta = \Theta_1 \i + \Theta_2 \j + \Theta_3 \k$. The $\Theta_i$ are defined only locally, but the so-called fundamental or Kraines real-valued  4-form $\Theta\wedge\bar{\Theta}$ is defined globally. The wedge stands here for the usual exterior product of quaternion-valued differential forms, whereas the bar stands for quaternionic conjugation; as $\Theta \in \mbox{Im\,}\mathbb{H}$, we have $\bar{\Theta} = -\Theta$. For the manifold to be quaternionic K\"ahler, the fundamental form has to be closed:
\begin{equation}
d(\Theta\wedge\bar{\Theta}) = 0 
\label{ff-closure}
\end{equation}
For $n>1$ this implies the existence of a locally defined $\mbox{Im\,}\mathbb{H}$-valued 1-form $\omega$ on ${\cal M}$ such that
\begin{equation}
d\Theta + \omega\wedge\Theta - \Theta\wedge\omega = 0 
\label{str-eq}
\end{equation}
This 1-form is just the $Sp(1)$ part of the Riemannian connection. All quaternionic K\"ahler manifolds of dimension greater than 4 are Einstein; this implies that the $Sp(1)$ part of the Riemannian curvature 2-form is proportional to $\Theta$:
\begin{equation}
d\omega + \omega\wedge\omega = s\,\Theta 
\label{Einstein-eq}
\end{equation}
The proportionality constant $s$ is the constant scalar curvature scaled by a dimension-dependent positive numerical factor.

For $n=1$ the closure condition \eqref{ff-closure} carries no meaning, since in four dimensions all 4-forms are closed trivially. Nevertheless, it can be shown that for Einstein self-dual spaces the Einstein property can be cast precisely in the form \eqref{Einstein-eq} if one takes $\omega$ to be the self-dual part of the spin connection and $\Theta$ to be a frame of $\wedge^2_-T^*{\cal M}$, the bundle of self-dual 2-forms on ${\cal M}$. Equation \eqref{str-eq} can be regarded in this context as a consistency condition for \eqref{Einstein-eq} rather than a solution for \eqref{ff-closure}.

\subsection*{Swann bundles}

The quaternionic K\"ahler structure of a manifold ${\cal M}$ can be canonically encoded into the hyperk\"ahler structure of a space with one extra quaternionic dimension -- its associated Swann bundle or hyperk\"ahler cone ${\cal U}({\cal M})$ \cite{MR1096180}. From a purely differential geometric point of view ${\cal U}({\cal M})$ is just ${\cal M}\times\mathbb{H}^{\times}$. Let $q = q_0 + q_1\i + q_2\j + q_3\k$ be the additional $\mathbb{H}$-valued coordinate. The hyperk\"ahler metric is then
\begin{equation}
G = s |q|^2 g + |dq - q\omega|^2
\end{equation}
with the $\mbox{Im\,}\mathbb{H}$-valued K\"ahler form
\begin{equation}
\Omega = s\pt q \Theta \bar{q} + (dq - q\omega) \!\wedge\! (\overline{dq - q\omega})
\end{equation}
A straightforward computation gives
\begin{equation}
d\Omega = d[\pt q (s\pt\Theta - d\omega - \omega\wedge\omega) \bar{q}\, ]
\end{equation}
which indeed manifestly vanishes when the Einstein property \eqref{Einstein-eq} holds.

There is a natural action of the group $\mathbb{H}^{\times}$ on ${\cal U}({\cal M})$ induced by the left quaternionic multiplication in the fiber, $q \mapsto uq$ with $u \in \mathbb{H}^{\times}$. If $u \in Sp(1)\subset \mathbb{H}^{\times}$, the metric remains invariant, otherwise it transforms conformally. Consider the left-invariant $\mathbb{H}$-valued 1-form
\begin{equation}
d\bar{q}\,\bar{q}^{-1\!}  = \sigma_0 + \sigma_1 \i + \sigma_2\pt \j + \sigma_3 \k
\end{equation}
Its components are given explicitly by 
{\allowdisplaybreaks
\begin{align}
\sigma_0 & = \frac{q_0dq_0+q_1dq_1+q_2dq_2+q_3dq_3}{|q|^2} \nonumber \\
\sigma_1 & = \frac{q_1dq_0-q_0dq_1+q_2dq_3-q_3dq_2}{|q|^2} \nonumber \\
\sigma_2 & = \frac{q_2dq_0-q_0dq_2+q_3dq_1-q_1dq_3}{|q|^2} \nonumber \\
\sigma_3 & = \frac{q_3dq_0-q_0dq_3+q_1dq_2-q_2dq_1}{|q|^2} \label{left-inv}
\end{align} 
}and form a closed algebra under the action of the de Rham operator:
\begin{align}
d \sigma_0 & = 0 \\[4pt]
d\sigma_i & = \epsilon_{ijk\,} \sigma_j\wedge\sigma_k
\end{align}
The component 1-forms can be alternatively expressed in an Euler coordinate basis through the transformation
\begin{equation}
q = r\, e^{-\i \frac{\phi}{2}}  e^{-\j \frac{\theta}{2}}  e^{-\k \frac{\psi}{2}} \label{Euler-coords}
\end{equation}
This gives
{\allowdisplaybreaks
\begin{align}
\sigma_0 & = \frac{dr}{r} \nonumber \\
\sigma_1 & = \frac{1}{2}(\cos\theta\cos\psi\, d\phi - \sin\psi\, d\theta) \nonumber \\
\sigma_2 & = \frac{1}{2}(\cos\theta\sin\psi\, d\phi + \cos\psi\, d\theta) \nonumber \\
\sigma_3 & = \frac{1}{2}(d\psi - \sin\theta\, d\phi)
\end{align} 
}and one can recognize the familiar form of the left-invariant Cartan-Maurer forms for the group $SU(2) \simeq Sp(1)$. 

Let us now observe that, for any $q \in \mathbb{H}^{\times}$, we have
$ |dq-q\omega|^2 = |q|^2|q^{-1\!}dq-\omega|^2 = |q|^2|d\bar{q}\,\bar{q}^{-1\!} \!+ \omega|^2$.
In the last step we have made use of the fact that $\omega$ is purely imaginary and so $\bar{\omega}=-\omega$. Based on this observation, the metric can be re-written as follows
\begin{equation}
G = |q^2| [s\pt g +  (\vec{\sigma} + \vec{\omega})^2  + \sigma_0^2] \label{G-sigma-basis}
\end{equation}
Passing to the radial coordinate, we can further write this as
\begin{equation}
G = dr^2 + r^2 [s\pt g +  (\vec{\sigma} + \vec{\omega})^2] \label{G-HKC}
\end{equation}
which is a standard cone metric. This form of the metric appears also in \cite{Santillan:2007qj}. The hyperk\"ahler 2-forms can be expressed as well in this basis; one has
\begin{equation}
\Omega = q [s\pt\Theta + (\sigma_0 - \vec{\sigma} - \vec{\omega})\wedge(\sigma_0 + \vec{\sigma} + \vec{\omega})] \bar{q}
\end{equation}
We have stretched a bit the notation here: $\vec{\sigma}$ and $\vec{\omega}$ should be regarded not as $\mathbb{R}^3$ vector-valued but of course as $\mbox{Im\,}\mathbb{H}$-valued 1-forms. The corresponding expressions for the components of $\Omega$ in Euler coordinates are straightforward to obtain.

The generators of the $\mathbb{H}^{\times}$ action on the Swann bundle are not the dual vector fields of the left-invariant one forms \eqref{left-inv} but rather those dual to the components of the right-invariant $\mathbb{H}$-valued 1-form $dq\,q^{-1}$. Explicitly,
{\allowdisplaybreaks
\begin{align}
X_0 & = q_0\frac{\partial}{\partial q_0} + q_1\frac{\partial}{\partial q_1} + q_2\frac{\partial}{\partial q_2} + q_3\frac{\partial}{\partial q_3} \nonumber \\
X_1 & = q_1\frac{\partial}{\partial q_0} - q_0\frac{\partial}{\partial q_1}  + q_3\frac{\partial}{\partial q_2} - q_2\frac{\partial}{\partial q_3} \nonumber \\
X_2 & = q_2\frac{\partial}{\partial q_0} - q_0\frac{\partial}{\partial q_2}  + q_1\frac{\partial}{\partial q_3} - q_3\frac{\partial}{\partial q_1} \nonumber \\
X_3 & = q_3\frac{\partial}{\partial q_0} - q_0\frac{\partial}{\partial q_3}  + q_2\frac{\partial}{\partial q_1} - q_1\frac{\partial}{\partial q_2} 
\end{align} 
}They satisfy, as required, the commutation relations of the $\mathbb{H}^{\times}$ algebra
\begin{align}
[X_i,X_0] & = 0 \\[4pt]
[X_i,X_j] & = 2\pt \epsilon_{ijk\,}X_k 
\end{align}
and one can check directly that indeed, the $Sp(1)$ generators $X_i$ act isometrically and rotate the hyperk\"ahler structures into one another, whereas $X_0$ acts conformally: 
\begin{align}
{\cal L}_{X_i}G & = 0 & {\cal L}_{X_i}\Omega_j & = 2\pt \epsilon_{ijk\,}\Omega_k \\[4pt]
{\cal L}_{X_0}G & = 2\pt G & {\cal L}_{X_0}\Omega_j & = 2\pt \Omega_j
\end{align}
In fact, since $X_0 = r \partial_r$ and given the conic form of the Swann bundle metric, $X_0$ generates a conformal homothety with the tip of the cone as its fixed point.

The action of each vector field $X_i$ is Hamiltonian with respect to the corresponding hyperk\"ahler symplectic form $\Omega_i$. Remarkably, all three associated moment maps coincide; they are equal, up to an additive constant, to $|q|^2$. According to a lemma proved by Hitchin \cite{Hitchin:1986ea,MR1726926}, this then implies that $|q|^2$ also plays the role of K\"ahler potential for all three standard (and, in fact, for the whole 2-sphere's worth of) complex structures of the Swann bundle; accordingly, one refers to it as hyperk\"ahler potential.

\subsection*{Tri-Hamiltonian actions}

Suppose on ${\cal U}({\cal M})$ there exists a free action generated by a vector field $X$ that commutes with the $\mathbb{H}^{\times}$-action and is tri-Hamiltonian with respect to the hyperk\"ahler symplectic forms, that is ${\cal L}_X\Omega = 0$. Then this action descends to an action $X_H$ on ${\cal M}$ which rotates the components of the almost K\"ahler 2-form $\Theta$, in the sense that there exists an $\mbox{Im\,}\mathbb{H}$-valued function $R$ on ${\cal M}$ such that 
\begin{equation}
{\cal L}_{X_H}\Theta = R\pt \Theta-\Theta R \label{QK-rot}
\end{equation}
or, equivalently, in $\mathbb{R}^3$ vector form, ${\cal L}_{X_H}\vec{\Theta} = 2\pt \vec{R}\times\vec{\Theta}$. Observe that any such action preserves the fundamental 4-form $\Theta \wedge \bar{\Theta}$ of the quaternionic K\"ahler space \cite{MR872143}. Conversely, any action $X_H$ on ${\cal M}$ which rotates the components of the 2-form $\Theta$ in the above sense can be lifted canonically, by combining it to an action on the $\mathbb{H}^{\times}$ fibers, to a tri-Hamiltonian action $X$ on ${\cal U}({\cal M})$ with corresponding $\mbox{Im\,}\mathbb{H}$-valued moment map
\begin{equation}
\mu_X = q\pt( i_{X_H}\omega+R)\pt\bar{q} \label{momma}
\end{equation}

To prove the direct implication, notice that the result of contracting any vector field $X$ from the tangent bundle of ${\cal U}({\cal M})$ with the $\mbox{Im\,}\mathbb{H}$-valued hyperk\"ahler 2-form of ${\cal U}({\cal M})$ can be cast in the form
\begin{align}
i_X\Omega {} = {} & q\pt[\pt i_X(s\pt\Theta - d\omega - \omega\wedge\omega) + {\cal L}_X\omega]\pt\bar{q}  - d(q\pt i_X\omega\pt\bar{q}) \nonumber \\[2pt]
& + (i_Xdq)(\overline{dq - q\omega}) - (dq-q\omega)(i_Xd\bar{q})
\end{align}
Assume further that $X$ commutes with the $\mathbb{H}^{\times}$-action on ${\cal U}({\cal M})$. This implies that $X$ admits a canonical decomposition into horizontal and vertical components as follows: $X = X_H+X_V$, with the $q$-independent $X_H$ acting on the base ${\cal M}$ and $X_V$ acting on the $\mathbb{H}^{\times}$ fibers such that $i_{X_V}dq = qR$, where $R$ is some arbitrary, $q$-independent, $\mbox{Im\,}\mathbb{H}$-valued function on ${\cal M}$. Using this in the above formula together with the Einstein property \eqref{Einstein-eq}, we get
\begin{equation}
i_X\Omega = q({\cal L}_{X_H}\omega + dR + \omega R - R\pt \omega)\bar{q} - d[\pt q\pt( i_X\omega + R)\pt\bar{q}\pt]  \label{ixO}
\end{equation}
Using moreover the fact that $\Omega$ is closed, one derives immediately the Lie action of $X$ on $\Omega$:
\begin{equation}
{\cal L}_X\Omega = d[\pt q({\cal L}_{X_H}\omega + dR + \omega R - R\pt \omega)\bar{q}\pt ]
\label{LxO}
\end{equation}
Clearly, ${\cal L}_X\Omega = 0$ if and only if
\begin{equation}
dR + \omega R - R\pt \omega = - {\cal L}_{X_H}\omega \label{R-eq}
\end{equation}
By acting on this equation with the de Rham operator and resorting again to the property \eqref{Einstein-eq} one arrives at the equation \eqref{QK-rot}.

The converse implication is now straightforward. Consider a vector field $X_H$ from $T{\cal M}$ having the property \eqref{QK-rot} for some $\mbox{Im\,}\mathbb{H}$-valued function $R$. Assemble the vector field $X=X_H+X_V$ in $T{\cal U}({\cal M})$, with fiberwise-acting $X_V$ defined by $i_{X_V}dq = qR$. Evidently, equations \eqref{ixO} and \eqref{LxO} hold for $X$. Notice on the other hand that equation \eqref{QK-rot} is the integrability condition for \eqref{R-eq}, so it follows that this one holds as well. From \eqref{ixO} and \eqref{LxO} one infers then that $X$ is tri-Hamiltonian with respect to $\Omega$, with corresponding moment maps given by the expression \eqref{momma}. Observe, as an aside, that if we define $\mu_{X_H} = i_{X_H}\omega+R$, then  equations \eqref{R-eq} and \eqref{Einstein-eq} imply that
\begin{equation}
d\mu_{X_H} + \omega\mu_{X_H} - \mu_{X_H}\omega = - s\, i_{X_H}\Theta
\end{equation}
This can be regarded as a generalization of the hyperk\"ahler moment map equation to the quaternionic K\"ahler case. 

These results extend those of \cite{MR1096180} and agree with the findings of \cite{deWit:2001bk}. The authors of \cite{deWit:2001bk} argue moreover that any isometry of a quaternionic K\"ahler space satisfies the rotation property \eqref{QK-rot} for some $R$ and can be lifted to a tri-holomorphic isometry on the Swann bundle.

\section{Toric hyperk\"ahler manifolds}

The general local formulation of $4(n+1)$-dimensional hyperk\"ahler metrics with a free  action of $T^{n+1}$ preserving the hyperk\"ahler structure has been given by Gibbons and Hawking for the case of four dimensions \cite{Gibbons:1979zt} and generalized to arbitrary dimensions by Pedersen and Poon \cite{MR953820}, building on the work of Hitchin, Karlhede, Lindstr\"om and Ro\v{c}ek \cite{Hitchin:1986ea}. 

Consider the trivial $\mathbb{R}^{n+1}$ bundle over $\mathbb{R}^{n+1}\!\otimes\mathbb{R}^3$, with connection 1-form $A = (A_0, \cdots, A_n)$ and the Higgs fields $\phi = (\phi_0, \cdots, \phi_n)$; $A_{\K}$ are 1-forms  on $\mathbb{R}^{n+1}\!\otimes\mathbb{R}^3$ with values in $\mathbb{R}$ and $\phi_{\K} = (\phi_{\K0}, \cdots, \phi_{\K n})$ are defined on $\mathbb{R}^{n+1}\!\otimes\mathbb{R}^3$ with values in $\mathbb{R}^{n+1}$. The pair $(A,\phi)$ is assumed to satisfy the following linear system of PDEs -- the generalized abelian monopole equations
\begin{align}
dA_{\K} = \star_{\I}d\Phi_{\K\I} && \vec{\nabla}_{\I}\Phi_{\K\J} = \vec{\nabla}_{\J}\Phi_{\K\I}
\label{gen-monopole-eqs}
\end{align}
In the first equation summation over the repeated index $I$ is understood. We define the linear Hodge-like operators $\star_{\I} : \mathbb{R}^{n+1} \!\otimes T^*\mathbb{R}^3 \to \mathbb{R}^{n+1}\! \otimes \wedge^2 T^*\mathbb{R}^3$ by specifying their action on a basis of $\mathbb{R}^{n+1}\! \otimes T^*\mathbb{R}^3$:
\begin{equation}
\star_{\I}d\vec{r}^{\,\J} = d\vec{r}^{\,\I} \!\wedge d\vec{r}^{\,\J}
\label{star-op}
\end{equation}
where
\begin{equation}
(d\vec{r}^{\,\I} \!\wedge d\vec{r}^{\,\J})_k = \frac{1}{2}\epsilon_{kij\,}dx^{\I}_i \!\wedge dx^{\J}_j
\end{equation}
is the standard vector product in $\mathbb{R}^3$. The hyperk\"ahler metric takes then the form of a generalized Gibbons-Hawking Ansatz
\begin{equation}
G = \frac{1}{2} \Phi_{\I\J} d\vec{r}^{\,\I} \!\!\cdot\! d\vec{r}^{\,\J} + \frac{1}{2} \Phi^{\I\J}(d\psi_{\I} + A_{\I}) (d\psi_{\J} + A_{\J}) 
\label{GH-Ansatz}
\end{equation}
with $\Phi^{\I\J}$ denoting the inverse of $\Phi_{\I\J}$. Moreover, the corresponding hyperk\"ahler 2-forms are \cite{Papadopoulos:1994kj,Gauntlett:1997pk}
\begin{equation}
\vec{\Omega} = \Phi_{\I\J} d\vec{r}^{\,\I} \!\wedge d\vec{r}^{\,\J} - (d\psi_{\I} + A_{\I}) \wedge d\vec{r}^{\,\I} 
\end{equation}
Indeed, by resorting to the second monopole equation \eqref{gen-monopole-eqs} one can show that $d\vec{\Omega} = (\star_{\I}d\Phi_{\K\I}-dA_{\K}) \wedge d\vec{r}^{\,\K}$, and this obviously vanishes provided that the first monopole equation holds.

The metric $G$ has $n+1$ isometries generated by the vector fields $X^{\I} = \partial_{\psi_I}$. Contraction with the hyperk\"ahler forms yields
\begin{equation}
i_{X^I} \vec{\Omega} = - d \vec{r}^{\,\I}
\label{mommap-coord}
\end{equation}
Since the three components of $\vec{\Omega}$ are closed, it follows that the action of $X^{\I}$ is tri-Hamiltonian (and hence tri-holomorphic), with corresponding moment maps $\vec{r}^{\,\I}$. Formula \eqref{GH-Ansatz} is therefore sometimes referred to as the moment-map basis description of the metric.

\section{The generalized Legendre transform \\ construction}

The generalized Legendre transform approach of Lindstr\"om and Ro\v{c}ek \cite{Lindstrom:1983rt} gives a solution to the monopole equations in terms of a single real-valued polyharmonic function $F$ on $\mathbb{R}^{n+1}\otimes\mathbb{R}^3$.  As shown by Bielawski \cite{MR1704547}, this solution is most generic, providing a complete local description of hyperk\"ahler metrics with free tri-holomorphic toric actions of rank equal to the quaternionic dimension. The construction emerged originally in physics in relation to supersymmetric field theories, as the superspace equivalent of Hodge duality between 0-form and 2-form gauge fields in four dimensions. Subsequently it acquired a geometric interpretation within the framework of the twistor theory of hyperk\"ahler manifolds \cite{Hitchin:1986ea}.

The function $F$ is constructed by contour-integrating a meromorphic function of sections $\hat{\eta}^{\I}$ of the pulled-back ${\cal O}(2)$ bundles over $Z$,  the twistor space of the hyperk\"ahler manifold:
\begin{equation}
F = \frac{1}{2\pi i} \oint_{\Gamma} \frac{d\zeta}{\zeta} {\cal H}(\hat{\eta}^{\I}) \label{GLT-action}
\end{equation}
The sections $\hat{\eta}^{\I}$ are required to satisfy a reality condition with respect to the real structure induced on $Z$ by antipodal conjugation on the Riemann sphere: $\overline{\hat{\eta}^{\I}(\zeta)} = \hat{\eta}^{\I}(-1/\bar{\zeta})$. With a slightly unusual choice of local trivialization, they take the form
\begin{equation}
\hat{\eta}^{\I}(\zeta) = \frac{\bar{z}^{\I}}{\zeta} + x^{\I} - z^{\I}\zeta \label{O2-multiplet}
\end{equation}
with $x^{\I} \in \mathbb{R}$. The presence of the real structure allows one to choose integration contours $\Gamma$ that result in real-valued functions $F$. By construction, $F$ is a function of the moduli of the ${\cal O}(2)$ sections. These are related to the vector-valued coordinates of the previous section by a complex-linear transformation,
\begin{equation}
\vec{r}^{\,\I} = 2\mbox{\pt Im\,}z^{\I}\,\i - 2\mbox{\pt Re\,}z^{\I}\,\j + x^{\I}\,\k
\end{equation}
Thus defined, $F$ automatically satisfies a system of linear second order PDEs, namely
\begin{align}
\Delta_{\I\J}F = 0 && F_{x_i^Ix_j^J} = F_{x_i^Jx_j^I}
\label{polyharm}
\end{align}
where $\Delta_{\I\J} = \vec{\nabla}_{\I} \sprod \vec{\nabla}_{\J}$ is a Laplacian-like operator and $x^{\I}_i$ with $i=1,2,3$ are the components of $\vec{r}^{\,\I}$; the indices of $F$ denote derivatives. Functions with this property are termed polyharmonic.

Lindstr\"om and Ro\v{c}ek show that the Legendre transform of $F$ with respect to the real coordinates $x^{\I}$ gives a K\"ahler potential ${\cal K}$ of the hyperk\"ahler metric as well as a corresponding set of local holomorphic coordinates $z^{\I}, u_{\I}$. Specifically,
\begin{equation}
{\cal K}(z^{\I},\bar{z}^{\I},u_{\I},\bar{u}_{\I}) = F(z^{\I},\bar{z}^{\I},x^{\I}) - x^{\I}(u_{\I} + \bar{u}_{\I})
\end{equation}
where the $x^{\I}$ are determined by
\begin{equation}
\frac{\partial F}{\partial x^{\I}} = u_{\I} + \bar{u}_{\I}
\end{equation}
These holomorphic coordinates are moreover Darboux coordinates for the hyperk\"ahler symplectic form, {\it i.e.},
\begin{eqnarray}
\Omega_+ = \frac{1}{2}(\Omega_1+i\Omega_2) = du_{\I} \!\wedge dz^{\I}
\end{eqnarray}
Purely imaginary shifts of $u_{\I}$ leave the Legendre relations invariant. It follows that these are isometric transformations, and so, one identifies
\begin{equation}
\psi_{\I} = \mbox{Im\,} u_{\I} \label{psi-coord}
\end{equation}
By extracting the metric from the K\"ahler potential and comparing the result with equation \eqref{GH-Ansatz}, one further identifies
\begin{align}
\Phi_{\I\J} & = -\frac{1}{2}F_{x^Ix^J} \label{Phi-GLT} \\
A_{\K}\, & = -\mbox{\pt Im\pt} (F_{x^Kz^J}dz^{\J}) \label{A-GLT}
\end{align}
It is a simple exercise to verify that these expressions provide indeed a solution to the generalized abelian monopole equations \eqref{gen-monopole-eqs} as long as $F$ satisfies the polyharmonicity conditions \eqref{polyharm}.

\subsection*{Swann bundles and the generalized Legendre transform}

One may ask the question: When does a hyperk\"ahler space constructed by means of the generalized Legendre transform have a Swann bundle structure? In \cite{deWit:2001dj}, de Wit, Ro\v{c}ek and Vandoren show that this happens provided that the meromorphic function ${\cal H}(\hat{\eta}^{\I})$ from which $F$ is built satisfies the following two requirements: \\
\phantom{\quad} $1)$ has no explicit dependence on $\zeta$ except via the $\hat{\eta}^{\I}$, and \\
\phantom{\quad} $2)$ is made up of either terms homogeneous of degree 1 in the $\hat{\eta}^{\I}$ or of the form $\hat{\eta}^{\I} \ln \hat{\eta}^{\I}$ (no summation over the index $I$ implied).

Observe that we have at our disposal a natural action of $\mathbb{H}^{\times}$ on the moduli space of ${\cal O}(2)$ sections, namely the one comprising rigid rotations and simultaneous rescalings of the monopole position vectors $\vec{r}^{\,\I}$, with generators $\vec{L} = -\vec{r}^{\,\I} \times \vec{\nabla}_{\I}$ and $L_0 = \vec{r}^{\,\I} \sprod \vec{\nabla}_{\I}$ (summation over the index $I$ implied). We will presently show that, provided that the above two conditions are satisfied, $\vec{L}$ and $L_0$ induce an $\mathbb{H}^{\times}$ action on the hyperk\"ahler space. 

The the two conditions on the function ${\cal H}(\hat{\eta}^{\I})$ translate into the following two linear differential equations for $F$,
\begin{equation}
L_3(F) = 0 \hspace{30pt} \mbox{and} \hspace{30pt} L_0(F) = 2F
\end{equation}
respectively \cite{Ionas:2008gh}. Based on these as well as on the polyharmonicity conditions \eqref{polyharm}, one can show that the Lie action of $\vec{L}$ and $L_0$ on the generalized Gibbons-Hawking Higgs fields and connection 1-forms yields
\begin{align}
{\cal L}_{\vec{L}\,}\Phi_{\I\J} & = 0 & {\cal L}_{\vec{L}\,}A_{\K} & = d(\, i_{\vec{L}}A_{\K} - \Phi_{\K\J}\vec{r}^{\,\J}) \label{rot-inv} \\[4pt]
{\cal L}_{L_0}\Phi_{\I\J} & = - \Phi_{\I\J} & {\cal L}_{L_0}A_{\K} & = 0 
\end{align}
Thus, under the action of $\vec{L}$ the Higgs fields remain invariant while the connections $A_{\K}$ get shifted with total derivatives; nevertheless, their gauge equivalence classes remain invariant, and this is all that matters, since any total derivative shifts can be absorbed into redefinitions of the $\psi_{\K}$. One concludes that $\vec{L}$ are Killing vectors and, similarly, that $L_0$ is a conformal Killing vector for the hyperk\"ahler metric \eqref{GH-Ansatz}. What this means is that the metric is invariant at rigid rotations of the monopole configuration and transforms conformally at simultaneous rescalings of the monopole position vectors.

\section{$\mathbb{H}^{\times}$-reduction of toric Swann bundles} \label{SEC:H*-red}

\subsection*{Collective vs. individual degrees of freedom and coordinates}

A $4(n+1)$-dimensional Swann bundle ${\cal U}({\cal M})$  with a free action of the torus $T^{n+1}$ that commutes with $\mathbb{H}^{\times}$ admits, as a toric hyperk\"ahler space, a description {\it \`{a} la} Pedersen and Poon. Our strategy is to match this description with the Swann picture in order to extract information about the geometry of the quaternionic K\"ahler $4n$-manifold ${\cal M}$. This is equivalent to putting into effect a quaternionic reduction program.

The tri-holomorphic toric isometries generated by the vector fields $X^{\I} = \partial_{\psi_I}$ commute with the $\mathbb{H}^{\times}$ action generators $\vec{L}$ and $L_0$ and so, according to the discussion from the last subsection of section \ref{SEC:QK-&-SB}, they descend to quaternionic isometries on the underlying quaternionic K\"ahler manifold. The corresponding hyperk\"ahler moment maps $\vec{r}^{\,\I}$ are $\mathbb{H}^{\times}$-invariant, meaning that under the action of some $q \in \mathbb{H}^{\times}$ they transform as follows:
\begin{equation}
\vec{r}^{\,\I} \hspace{5pt} \stackrel{q \in \mathbb{H}^{\times}}{\longrightarrow} \hspace{5pt} \vec{r}^{\,\I \prime} = q\,\vec{r}^{\,\I}\bar{q} \label{H*-equiv}
\end{equation}
for all $I$. This substantiates the intuitive picture of collective transformations of monopole configurations that emerged earlier. In writing the equation on the right hand side of \eqref{H*-equiv} we have committed a little abuse of notation which, unrepentantly, we turn into a policy: throughout this paper, for expediency reasons, we will occasionally blur the line between $\mathbb{R}^3$ vectors and imaginary quaternions whenever the context allows for an unequivocal interpretation. For example, in $\vec{r}^{\,\I \prime} = q\,\vec{r}^{\,\I}\bar{q} $, despite the notation's suggestion to the contrary, we clearly regard $\vec{r}^{\,\I}$ and $\vec{r}^{\,\I\prime}$ as imaginary quaternions and not as $\mathbb{R}^3$ vectors.

We can use the equivariance property to coordinatize the quaternionic K\"ahler manifold. This is essentially a local procedure, and locally ${\cal U}({\cal M})$ is isomorphic to the direct product ${\cal M}\times\mathbb{H}^{\times}$. Since the toric isometries $X^{\I}$ descend from ${\cal U}({\cal M})$ to ${\cal M}$, one can use the corresponding $\psi_{\I}$ coordinates as coordinates on ${\cal M}$. But ${\cal M}$ is $4n$-dimensional, so these do not suffice; we still need to define $3n-1$ more coordinates. To do that, let us observe that if we associate to a point $p \in {\cal M}$ a fixed configuration of the monopoles, then the collective transformations of this configuration (that is, rigid rotations and simultaneous scalings) are encoded in the $\mathbb{H}^{\times}$ fiber on top of $p$. Thus, we can think of the points of ${\cal M}$ as parametrizing inequivalent configurations of $n+1$ monopoles. Such a configuration has $3n-1$ degrees of freedom: $3(n+1)$ individual minus 4 collective, exactly the number of coordinates we need. To make things more concrete, we choose the reference configurations such that
{\allowdisplaybreaks
\begin{align}
\vec{r}^{\,\0} & = q\pt \k\pt \bar{q} \label{r_0} \\[4pt]
\vec{r}^{\,\1} & = q\pt(\rho^{\1}\j + \eta^{\1}\k)\pt\bar{q} \label{r_1} \\[4pt]
\vec{r}^{\,\I} & = q\pt(\chi^{\I}\i + \rho^{\I}\j + \eta^{\I}\k)\pt\bar{q} \label{r_I}
% \hspace{12pt} \mbox{for $I \geq 2$}
\end{align}
}for $I \geq 2$. This choice is arbitrary, without being restrictive. We will sometimes refer to this set of equations compactly as
\begin{equation}
\vec{r}^{\,\I} = q\, \vec{\chi}^{\,\I}\bar{q} \label{r=qchibq}
\end{equation}
with $\vec{\chi}^{\,\I} = \chi^{\I}\i + \rho^{\I}\j + \eta^{\I}\k$ such that $\chi^{\1}=\chi^{\0}=\rho^{\0}=0$ and $\eta^{\0}=1$. The remaining components of the vectors $\vec{\chi}^{\,\I}$ provide, together with the toric angles $\psi_{\I}$, a concrete set of local coordinates on ${\cal M}$.

\subsection*{The $q$-dependence of Higgs fields}

Our next step will be to re-express the two monopole equations \eqref{gen-monopole-eqs}, given in the ``homogeneous" coordinate basis $\vec{r}^{\,\I}$, in terms of the equivalent basis consisting of  the ``inhomogeneous" base coordinates $\vec{\chi}^{\,\I}$ and fiber coordinates $q$. In effect, this amounts to performing a reduction of the monopole equations along the Swann fiber bundle structure.

Let us observe that the first equation \eqref{rot-inv}, asserting the  invariance of the Higgs fields $\Phi_{\I\J}$ under collective rotations of the monopole configuration, implies that these must depend on the $q$-coordinates only through the norm $|q|^2$. The equation underneath shows moreover that $\Phi_{\I\J}$ scales at an overall scaling of the monopole position vectors with weight $-1$. This means that the Higgs fields must be of the form
\begin{equation}
\Phi_{\I\J} = \frac{U_{\I\J}}{|q|^2}  
\label{Higgs-Coulomb}
\end{equation}
with $U_{\I\J}$ a field defined entirely on ${\cal M}$, with no $q$-dependence.

\subsection*{Reduction of the second monopole equation}

We begin by examining the second monopole equation \eqref{gen-monopole-eqs}. In view of equation  \eqref{Higgs-Coulomb}, the gradients of the Higgs field components can be shown to take the following form in the inhomogeneous coordinate basis:
{\allowdisplaybreaks
\begin{align}
\vec{\nabla}_{\0}\Phi_{\K\J} & = - \frac{q ( A_{\K\J} \i + B_{\K\J} \j + C_{\K\J} \k) \bar{q} }{|q|^6} \nonumber \\
\vec{\nabla}_{\1}\Phi_{\K\J} & = \frac{q ( D_{\K\J} \i + \partial_{\rho^1}U_{\K\J} \j +\partial_{\eta^1}U_{\K\J} \k) \bar{q} }{|q|^6} \nonumber \\
\vec{\nabla}_{\I}\Phi_{\K\J} & = \frac{q ( \partial_{\chi^{I}}U_{\K\J} \i + \partial_{\rho^{I}}U_{\K\J} \j +\partial_{\eta^{I}}U_{\K\J} \k) \bar{q} }{|q|^6} 
\end{align} 
}for $I \geq 2$ and all values of $K,J$, where
{\allowdisplaybreaks
\begin{align}
A_{\K\J} & = (\eta^{\I}\partial_{\chi^{I}}-\chi^{\I}\partial_{\eta^I})U_{\K\J} + \frac{\eta^{\1}}{\rho^{\1}}(\chi^{\I}\partial_{\rho^{I}}-\rho^{\I}\partial_{\chi^I})U_{\K\J} \nonumber  \\
B_{\K\J} & = (\eta^{\I}\partial_{\rho^{I}}-\rho^{\I}\partial_{\eta^I})U_{\K\J} \nonumber \\[5pt]
C_{\K\J} & = (1+\chi^{\I}\partial_{\chi^{I}}+\rho^{\I}\partial_{\rho^I}+\eta^{\I}\partial_{\eta^I})U_{\K\J}  \nonumber \\[0pt]
D_{\K\J} & = \frac{1}{\rho^{\1}}(\chi^{\I}\partial_{\rho^{I}}-\rho^{\I}\partial_{\chi^I})U_{\K\J}
\end{align}
}To simplify appearances, here and throughout these notes we use the following convention: the indices $I,J,K, \cdots$ run over the maximum allowed range for which the expressions in which they appear make sense and are well-defined. For example, in $\eta^{\I}\partial_{\chi^I}$ the index $I$ runs from $2$ to $n$, whereas in $\eta^{\I}\partial_{\rho^I}$ it runs from $1$ to $n$. 

Using these expressions, it is straightforward to read off the consequences of the second monopole equation, $\vec{\nabla}_{\I}\Phi_{\K\J}  = \vec{\nabla}_{\J}\Phi_{\K\I}$. We get
\begin{empheq}[box=\fbox]{align}
\partial_{\chi^I}U_{\K\J} & = \partial_{\chi^J}U_{\K\I} \label{dU-sym1} \\[4pt]
\partial_{\rho^I}U_{\K\J} & = \partial_{\rho^J}U_{\K\I} \\[4pt]
\partial_{\eta^I}U_{\K\J} & = \partial_{\eta^J}U_{\K\I} \label{dU-sym3}
\end{empheq}
for all values of $K$ and the maximum allowed ranges for $I$ and $J$ (which are $I,J = \overline{2,n}$ in the first equation and $I,J = \overline{1,n}$ in the second and third). We also obtain
\begin{align}
& \partial_{\chi^I}U_{\K\0} = - A_{\K\I} && \partial_{\rho^I}U_{\K\0} = - B_{\K\I}  && \partial_{\eta^I}U_{\K\0} = - C_{\K\I} \label{coco} \\[4pt]  
& \partial_{\chi^I}U_{\K\1} = \phantom{+} D_{\K\I} \label{chanel}
\end{align}
as well as
\begin{equation}
A_{\K\1}+D_{\K\0} = 0
\end{equation}
again, for all values of $K$ and the maximum allowed ranges for $I$. 

One can show, with some effort, that this last equation is a consequence of the others. On the other hand, the equations involving $A_{\K\I}$, $B_{\K\I}$ and $D_{\K\I}$ in \eqref{coco}--\eqref{chanel} can be manipulated with the help of \eqref{dU-sym1}--\eqref{dU-sym3} into the following forms
\begin{empheq}[box=\fbox]{align}
\partial_{\chi^I}(U_{\K\J}\rho^{\J}) & = \partial_{\rho^I}(U_{\K\J}\chi^{\J}) \label{dB-sym1} \\[4pt]
\partial_{\rho^I}(U_{\K\J}\eta^{\J}) & = \partial_{\eta^I}(U_{\K\J}\rho^{\J}) \\[4pt]
\partial_{\eta^I}(U_{\K\J}\chi^{\J}) & = \partial_{\chi^I}(U_{\K\J}\eta^{\J}) \label{dB-sym3}
\end{empheq}
The remaining equation -- the one involving $C_{\K\I}$ -- can be similarly recast as
\begin{equation}
\boxed{ \chi^{\J}\partial_{\chi^J}U_{\K\I} + \partial_{\rho^I}(U_{\K\J}\rho^{\J}) + \partial_{\eta^I}(U_{\K\J}\eta^{\J}) = U_{\K\I} } \label{partial-tr}
\end{equation}
for all $I = \overline{1,n}$. Summing up, the second monopole equation \eqref{gen-monopole-eqs} on ${\cal U}({\cal M})$ reduces to equations \eqref{dU-sym1}--\eqref{dU-sym3}, \eqref{dB-sym1}--\eqref{dB-sym3} and \eqref{partial-tr} on ${\cal M}$.

\subsection*{Reduction of the first monopole equation}

To reduce the first monopole equation \eqref{gen-monopole-eqs} we need to evaluate the action of the $\star_{\I}$ operators on the inhomogeneous $dq$, $d\vec{\chi}^{\,\J}$ basis of $\mathbb{R}^{n+1}\otimes T^*\mathbb{R}^3$; since by definition the they act linearly, knowledge of their action on a basis is enough to completely determine them. This can be achieved by a direct calculation, but the resulting expressions turn out to be quite sizable and messy. Fortunately, we can benefit from a slight change of perspective. A significant simplification occurs if one uses instead of $dq = (dq_0,dq_i)$ the equivalent basis provided by the left $\mathbb{H}^{\times}$-invariant 1-forms $\sigma_0, \sigma_i$, defined in \eqref{left-inv}. In this new basis, equation \eqref{star-op} reads
{\allowdisplaybreaks
\begin{align}
\star_{\I}\sigma_0 = \phantom{+} \frac{|q|^2}{2\ } [& d\chi^{\I}\!\wedge\sigma_1 - 2\chi^{\I}(\sigma_1\wedge\sigma_0 - \sigma_2\wedge\sigma_3) \,+ \nonumber \\ 
& d\rho^{\I}\!\wedge\sigma_2 - 2\rho^{\I}(\sigma_2\wedge\sigma_0 - \sigma_3\wedge\sigma_1) + 4\eta^{\I}\sigma_1\wedge\sigma_2] \nonumber \\
\star_{\I}\sigma_1 = -\frac{|q|^2}{2\ } [& d\chi^{\I}\!\wedge\sigma_0 + 2 \rho^{\I}(\sigma_3\wedge\sigma_0 \!- \sigma_1\wedge\sigma_2) + d\eta^{\I}\!\wedge\sigma_2 - 4 \eta^{\I}\sigma_2\wedge\sigma_0] \nonumber \\[-1pt]
\star_{\I}\sigma_2 = -\frac{|q|^2}{2\ } [& d\rho^{\I}\!\wedge\sigma_0 - 2 \chi^{\I}(\sigma_3\wedge\sigma_0 \!- \sigma_1\wedge\sigma_2) - d\eta^{\I}\!\wedge\sigma_1 + 4 \eta^{\I}\sigma_1\wedge\sigma_0] \nonumber \\[-1pt]
\star_{\I}\sigma_3 = -\frac{|q|^2}{2\ } [& d\eta^{\I}\!\wedge\sigma_0 + 2 \chi^{\I}(\sigma_2\wedge\sigma_0 \!- \sigma_3\wedge\sigma_1) + d\rho^{\I}\!\wedge\sigma_1 - 4 \rho^{\I}\sigma_1\wedge\sigma_0 \nonumber \\[0pt]
& + \alpha^{\I}_3\!\wedge d\rho^{\1}/ \rho^{\1}  -  \alpha^{\I}_2\!\wedge d\eta^{\1}/ \rho^{\1} ]
\end{align} 
}and
{\allowdisplaybreaks
\begin{align}
\star_{\I}d\chi^{\J} = -|q|^2 [& \alpha^{\I}_3\!\wedge(d\rho^{\J} \!- \rho^{\J}d\rho^{\1}/\rho^{\1}) - \alpha^{\I}_2\!\wedge(d\eta^{\J} \!- \rho^{\J}d\eta^{\1}/\rho^{\1}) \nonumber \\[2pt] 
& + (\beta^{\I}_{11} - \beta^{\I}_{33})\chi^{\J} ] \nonumber \\[2pt] 
\star_{\I}d\rho^{\J} = \phantom{+} |q|^2 [& \alpha^{\I}_3\!\wedge(d\chi^{\J} \!- \chi^{\J}d\rho^{\1}/\rho^{\1}) \!+ \alpha^{\I}_2\!\wedge\chi^{\J}d\eta^{\1}/\rho^{\1} - \alpha^{\I}_1\!\wedge d\eta^{\J} \nonumber \\[2pt] 
& - (\beta^{\I}_{12}+\beta^{\I}_{21})\chi^{\J} - (\beta^{\I}_{22}-\beta^{\I}_{33})\rho^{\J} ] \nonumber \\[2pt] 
\star_{\I}d\eta^{\J} = -|q|^2 [& \alpha^{\I}_2\!\wedge d\chi^{\J} - \alpha^{\I}_1\!\wedge d\rho^{\J} \nonumber \\[2pt] 
& + (\beta^{\I}_{13}+\beta^{\I}_{31})\chi^{\J} \!+ (\beta^{\I}_{23}+\beta^{\I}_{32})\rho^{\J} ]
\end{align} 
}where, for brevity, we have introduced the notations
\begin{align}
\alpha^{\I}_i & = d\chi^{\I}_i + \chi^{\I}_i \sigma_0 + \epsilon_{ijk\,}\chi^{\I}_j\sigma_k \\[4pt]
\beta^{\I}_{ij} & = d\chi^{\I}_i \!\wedge\sigma_j - \chi^{\I}_i (2\sigma_j\wedge\sigma_0 + \epsilon_{jkl\,}\sigma_k\wedge\sigma_l)
\end{align}
Based on these relations and the linearity of the $\star$-operators we can evaluate their action on $d\Phi_{\K\I}$. An extremely long and laborious calculation in which we make extensive use of the equations that follow from reducing the second monopole equation \eqref{gen-monopole-eqs} is rewarded by the remarkably simple result
\begin{equation}
\star_{\I}d\Phi_{\K\I} = F_{\K} + d(\, \vec{\sigma}\pt \sprod \vec{B}_{\K}) 
\label{star-dPhi}
\end{equation}
where
\begin{equation}
\vec{B}_{\K} = 2\pt U_{\K\I}\vec{\chi}^{\,\I} \label{B-def}
\end{equation}
and
{\allowdisplaybreaks
\begin{empheq}[box=\fbox]{align}
F_{\K} & = \frac{1}{2}(\partial_{\chi^I}U_{\K\J}+\partial_{\chi^J}U_{\K\I})\, d\eta^{\I} \!\wedge d\rho^{\J} \nonumber \\[-1pt]
& + \,\frac{1}{2}(\partial_{\rho^I}U_{\K\J}+\partial_{\rho^J}U_{\K\I})\, d\chi^{\I} \!\wedge d\eta^{\J} \nonumber \\[-1pt]
& + \,\frac{1}{2}(\partial_{\eta^I}U_{\K\J}+\partial_{\eta^J}U_{\K\I})\, d\rho^{\I} \!\wedge d\chi^{\J} \nonumber \\[-1pt]
& + \,\frac{1}{2}(\chi^{\I}\partial_{\rho^I}U_{\K\J}-\rho^{\I}\partial_{\chi^I}U_{\K\J})\frac{d\eta^{\1} \!\wedge d\rho^{\J}+d\eta^{\J} \!\wedge d\rho^{\1}}{\rho^{\1}}
\end{empheq} 
}Summation over the index $I$ in the l.h.s. of \eqref{star-dPhi} is of course understood. Incidentally, let us mention that unlike the rest of the equation, the above expression for $F_{\K}$ is obtained directly, without any reference to the second monopole equation or its consequences.

The first monopole equation, $dA_{\K} = \star_{\I}d\Phi_{\K\I}$, readily implies then that the connection 1-forms $A_{\K}$ must be, up to inconsequential exact terms, of the form
\begin{equation}
A_{\K} = C_{\K} + \vec{\sigma}\pt \sprod \vec{B}_{\K} \label{ABC}
\end{equation}
with $C_{\K}$ 1-forms on ${\cal M}$ satisfying the reduced Bogomol'nyi equation
\begin{equation}
\boxed{ dC_{\K} = F_{\K} } \label{Bogo}
\end{equation}

\subsection*{The solution of the $\vec{B}_{\K}$ field equations}

In this and the following two subsections we will concern ourselves with solving the system of equations for $U_{\I\J}$ that results from reducing the second monopole equation. We will show that its solutions are determined by a single function $\V$ of the inhomogeneous coordinates $\vec{\chi}^{\,\I}$ satisfying a set of linear second-order differential constraints.

For that, it is useful to investigate first the properties of the vectors $\vec{B}_{\K}$ that emerged in connection to the first monopole equation. Let $\vec{B}_{\K} = {(B_{\K})}_1 \i + {(B_{\K})}_2\pt \j + {(B_{\K})}_3 \k$. The symmetry property of $U_{\I\J}$ implies that the components of $\vec{B}_{\K}$ must satisfy
\begin{align}
\chi^{\I}{(B_{\I})}_2 & = \rho^{\I}{(B_{\I})}_1 \label{U-sym1} \\[4pt]
\rho^{\I}{(B_{\I})}_3 & = \eta^{\I}{(B_{\I})}_2 \\[4pt]
\eta^{\I}{(B_{\I})}_1 & = \chi^{\I}{(B_{\I})}_3 \label{U-sym3}
\end{align}
Equations \eqref{dB-sym1}--\eqref{dB-sym3} can be rephrased in terms of $\vec{B}_{\K}$ as follows:
\begin{align}
\partial_{\chi^I}{(B_{\K})}_2 & = \partial_{\rho^I}{(B_{\K})}_1 \label{curlxB1} \\[4pt]
\partial_{\rho^I}{(B_{\K})}_3 & = \partial_{\eta^I}{(B_{\K})}_2 \label{curlxB2} \\[4pt]
\partial_{\eta^I}{(B_{\K})}_1 & = \partial_{\chi^I}{(B_{\K})}_3 \label{curlxB3}
\end{align}
Furthermore, the equations \eqref{dU-sym1}--\eqref{dU-sym3} imply that 
\begin{align}
\partial_{\chi^I}{(B_{\J})}_k & = \partial_{\chi^J}{(B_{\I})}_k \label{dchiB} \\[4pt]
\partial_{\rho^I}{(B_{\J})}_k & = \partial_{\rho^J}{(B_{\I})}_k \label{drhoB} \\[4pt]
\partial_{\eta^I}{(B_{\J})}_k & = \partial_{\eta^J}{(B_{\I})}_k \label{detaB}
\end{align}
for $k=1,2,3$. Eventually, multiplying successively equation \eqref{partial-tr} with $\chi^{\K}$, $\rho^{\K}$ and $\eta^{\K}$, summing up over the index $K$ and then making use of the equations \eqref{U-sym1}--\eqref{U-sym3}, we obtain, for all $k$,
\begin{equation}
\chi^{\J}\partial_{\chi^J}{(B_{\I})}_k + \rho^{\J}\partial_{\rho^I}{(B_{\J})}_k + \eta^{\J}\partial_{\eta^I}{(B_{\J})}_k = 0 \label{trace-part}
\end{equation}
In line with our previously declared convention, in equations \eqref{curlxB1} through \eqref{trace-part} the indices $I,J,K$ run over the maximum ranges for which these equations make sense. Thus, in \eqref{curlxB1} and \eqref{curlxB3} $I = \overline{2,n}$ while in \eqref{curlxB2} $I = \overline{1,n}$; in all of them $K = \overline{0,n}$. In \eqref{dchiB} $I,J = \overline{2,n}$ whereas in \eqref{drhoB} and \eqref{detaB} $I,J = \overline{1,n}$. Finally, in \eqref{trace-part} $I = \overline{1,n}$.

By resorting to the equations \eqref{curlxB1} through \eqref{detaB}, the equations \eqref{trace-part} (apart from the one with $k=1$, $I=1$) can be re-cast in the following form
\begin{align}
{(B_{\I})}_1 & = \partial_{\chi^I}\V \label{B1=d1V} \\[4pt]
{(B_{\I})}_2 & = \partial_{\rho^I}\V \\[4pt]
{(B_{\I})}_3 & = \partial_{\eta^I}\V \label{B3=d3V}
\end{align}
where
\begin{equation}
\V = \vec{B}_{\K}\sprod\pt\vec{\chi}^{\,\K} = 2\pt U_{\J\K} \vec{\chi}^{\,\J} \!\sprod \vec{\chi}^{\,\K}
\label{V-def}
\end{equation}
The indices $I$ can take the maximum allowed values, but this does not cover their entire range. In fact, we can show that all the $(B_{\I})_k$ can be expressed in terms of the function $V$. The remaining components can be determined from
{\allowdisplaybreaks
\begin{align}
\rho^{\1} {(B_{\1})}_1 & = (\chi^{\I}\partial_{\rho^I} -\rho^{\I}\partial_{\chi^I})\V \label{B1_1} \\[4pt]
{(B_{\0})}_1+ \eta^{\1} {(B_{\1})}_1 & = (\chi^{\I}\partial_{\eta^I} -\eta^{\I}\partial_{\chi^I})\V \\[4pt]
{(B_{\0})}_2 & = (\rho^{\I}\partial_{\eta^I} -\eta^{\I}\partial_{\rho^I})\V \\[4pt]
{(B_{\0})}_3 & = (1-\chi^{\I}\partial_{\chi^I} -\rho^{\I}\partial_{\rho^I} -\eta^{\I}\partial_{\eta^I})\V \label{B0_3}
\end{align} 
}The first three equations follow from \eqref{U-sym1}--\eqref{U-sym3} and the last one from \eqref{V-def}.

\subsection*{The solution of the reduced Higgs field equations}

The next step is to show that all the $U_{\I\J}$ can be expressed in terms of $\vec{B}_{\K}$ and its derivatives. By using the symmetry property \eqref{dU-sym1} in equation \eqref{partial-tr} with $I = I' \geq 2$ and $K$ unrestricted, we get
\begin{equation}
U_{\K\I'} = \frac{1}{4} [\partial_{\chi^{\I'}}{(B_{\K})}_1 + \partial_{\rho^{\I'}}{(B_{\K})}_2 + \partial_{\eta^{\I'}}{(B_{\K})}_3] \label{UB1}
\end{equation}
A prime over an index indicates that the respective index takes values from $2$ to $n$. Furthermore, the definition \eqref{B-def} can be turned around to yield
\begin{align}
2 U_{\K\1}\rho^{\1} & = {(B_{\K})}_2 - 2 U_{\K\I'} \rho^{\I'}  \\[3pt]
2 U_{\K\0} + 2 U_{\K\1}\eta^{\1} & = {(B_{\K})}_3 - 2 U_{\K\I'} \eta^{\I'} \label{UB3}
\end{align}
for all $K$. These two relations can be straightforwardly solved for $U_{\K\1}$ and $U_{\K\0}$ through a linear transformation.

Since, as shown above, all the $\vec{B}_{\K}$ are determined by the function $\V$, it follows that all the $U_{\I\J}$ are in turn determined entirely by $\V$ as well. We can in fact write the dependence relations explicitly. We have found that the formulas take a somehow simpler form when expressed in an alternate frame.  Letting $\tilde{U}_{\I\J} = e_{\I}{}^{\K}e_{\J}{}^{\scriptscriptstyle L}\pt U_{\K{\scriptscriptstyle L}}$, with $e_{\I}{}^{\J}$ equal to $\delta_{\I}{}^{\J}$ if $I \geq 2$, to $\rho^{\1}\pt \delta_1{}^{\J}$ if $I=1$ and to $\eta^{\0}\pt \delta_0{}^{\J} + \eta^{\1}\pt \delta_1{}^{\J}$ if $I=0$,
%\begin{displaymath}
%e_{\I}{}^{\J} = 
%\begin{cases} 
%\delta_{\I}{}^{\J} & \mbox{if\ \ } I \geq 2  \\ 
%\rho^{\1}\pt \delta_1{}^{\J} & \mbox{if\ \ }  I = 1 \\ 
%\eta^{\1}\pt \delta_1{}^{\J} + \eta^{\0}\pt \delta_0{}^{\J} & \mbox{if\ \ } I = 0
%\end{cases}
%\end{displaymath}
we have
{\allowdisplaybreaks
\begin{align}
\tilde{U}_{\I'\!\J'} & = U_{\I'\!\J'} \nonumber \\[5pt]
\tilde{U}_{\I'\1} & = \frac{1}{2}\partial_{\rho^{I'}}\V - U_{\I'\!\J'}\rho^{\J'} \nonumber \\[0pt]
\tilde{U}_{\I'\0} & = \frac{1}{2}\partial_{\eta^{I'}}\V - U_{\I'\!\J'}\eta^{\J'} \nonumber \\[0pt]
\tilde{U}_{\1\1\,} & = \frac{1}{2}(\rho^{\1}\partial_{\rho^1} - \rho^{\I'}\partial_{\rho^{I'}})\V + \rho^{\I'}U_{\I'\!\J'}\rho^{\J'} \nonumber \\[0pt]
\tilde{U}_{\0\1\,} & = \frac{1}{2}(\rho^{\1}\partial_{\eta^1} - \eta^{\I'}\partial_{\rho^{I'}})\V + \eta^{\I'}U_{\I'\!\J'}\rho^{\J'} \nonumber \\[0pt]
\tilde{U}_{\0\0} & = \frac{1}{2}(1 - \rho^{\1}\partial_{\rho^1} - \chi^{\I'}\partial_{\chi^{I'}} - \rho^{\I'}\partial_{\rho^{I'}} - 2\pt \eta^{\I'}\partial_{\eta^{I'}})\V + \eta^{\I'}U_{\I'\!\J'}\eta^{\J'}
\label{UV}
\end{align} 
}with
\begin{equation}
U_{\I'\!\J'} = \frac{1}{4}(\partial_{\chi^{I'}}\partial_{\chi^{J'}}+\partial_{\rho^{I'}}\partial_{\rho^{J'}}+\partial_{\eta^{I'}}\partial_{\eta^{J'}})\V 
\end{equation} 
The explicit expressions for $U_{\I\J}$ in terms of $\V$ can be retrieved from these relations by reverting to the initial frame.

\subsection*{Differential constraints on $\V$}

The differential equations that the reduced Higgs fields $U_{\I\J}$ satisfy impose constraints on $\V$. To determine these, it is advantageous to begin by finding instead the constraints imposed on $\V$ by the $\vec{B}_{\K}$ differential equations \eqref{curlxB1} through \eqref{trace-part}.

The symmetry properties \eqref{dchiB}--\eqref{detaB} together with relations \eqref{B1=d1V}--\eqref{B3=d3V} yield the following constraints (the indices $I,J$ run as usual over the maximum allowed ranges):
\begin{empheq}[box=\fbox]{align}
\partial_{\chi^I}\partial_{\rho^J}\V & = \partial_{\chi^J}\partial_{\rho^I}\V \label{ddV1} \\[4pt]
\partial_{\rho^I}\partial_{\eta^J}\V & = \partial_{\rho^J}\partial_{\eta^I}\V \label{ddV2} \\[4pt]
\partial_{\eta^I}\partial_{\chi^J}\V & = \partial_{\eta^J}\partial_{\chi^I}\V \label{ddV3}
\end{empheq}
To derive the conditions imposed on $\V$ by the closure relations \eqref{curlxB1}--\eqref{curlxB3} we proceed in a roundabout way. Define
\begin{equation}
W_{ij} = 2\pt U_{\I\J}\chi_i^{\I}\chi_j^{\J} + \delta_{ij}\V
\end{equation}
By re-writing this as ${(B_{\J})}_i\chi_j^{\J} + \delta_{ij}\V$ and using the relations  \eqref{curlxB1}--\eqref{curlxB3}, the symmetry of $U_{\I\J}$ and, for $i=j$, the expressions \eqref{B1=d1V}--\eqref{B3=d3V}, one can show that they satisfy the equations
\begin{empheq}[box=\fbox]{align}
\partial_{\rho^I}W_{1j} & = \partial_{\chi^I}W_{2j} \label{dubya1} \\[4pt]
\partial_{\eta^I}W_{2j} & = \partial_{\rho^I}W_{3j} \label{dubya2} \\[4pt]
\partial_{\chi^I}W_{3j} & =  \partial_{\eta^I}W_{1j} \label{dubya3}
\end{empheq}
On the other hand, using the expressions for the components of $\vec{B}_{\K}$ in terms of $\V$, we get
\begin{align}
& W_{11} = \chi^{\I}\partial_{\chi^I}\V+\V && W_{12} = \chi^{\I}\partial_{\rho^I}\V && W_{13} = \chi^{\I}\partial_{\eta^I}\V \\[4pt] 
& W_{21} = \chi^{\I}\partial_{\rho^I}\V && W_{22} = \rho^{\I}\partial_{\rho^I}\V + \V && W_{23} = \rho^{\I}\partial_{\eta^I}\V \\[4pt] 
& W_{31} = \chi^{\I}\partial_{\eta^I}\V && W_{32} = \rho^{\I}\partial_{\eta^I}\V && W_{33} = 2\V - \chi^{\I}\partial_{\chi^I}\V - \rho^{\I}\partial_{\rho^I}\V
\end{align}
Equations \eqref{dubya1}--\eqref{dubya3} can therefore be interpreted as a set of constraints for $\V$. Not all of these constraints are independent of the previous ones, though. For example, we have
\begin{equation}
\partial_{\rho^{I}}W_{32} - \partial_{\eta^{I}}W_{22} = \rho^{\J}(\partial_{\rho^{I}}\partial_{\eta^{J}}\V - \partial_{\rho^{J}}\partial_{\eta^{I}}\V)
\end{equation}
which vanishes automatically by way of \eqref{ddV2}. A maximal independent subset of constraints is given by the equations \eqref{dubya1}--\eqref{dubya3} with $j=3$ together with equation \eqref{dubya1} with $j=2$.

Finally, the equations \eqref{trace-part} give no new constraints. Apart from the one with $k=1$, $I=1$, they have been used to derive the relations \eqref{B1=d1V}--\eqref{B3=d3V}. As for this one, after multiplication with an overall factor $\rho^{\1}$ it can be re-written as follows
\begin{equation}
\rho^{\I'}(\partial_{\chi^{\I'}}W_{33}-\partial_{\eta^{\I'}}W_{13}) + \rho^{\I'}(\partial_{\chi^{\I'}}W_{22}-\partial_{\rho^{\I'}}W_{12}) + \chi^{\I'}(\partial_{\eta^{\I'}}W_{23}-\partial_{\rho^{\I'}}W_{33}) = 0
\end{equation}
and this obviously holds identically in view of the previous constraints.

This exhausts the list of $\vec{B}_{\K}$ properties. One can furthermore show that the $U_{\I\J}$ properties \eqref{dU-sym1}--\eqref{dU-sym3} impose no new constraints on $\V$ either. In fact, one can prove that if we use the relations \eqref{UB1} through \eqref{UB3} to define $U_{\I\J}$ in terms of $\vec{B}_{\K}$, then the equations \eqref{dU-sym1}--\eqref{dU-sym3} with $K \geq 1$ follow as consequences of the properties of $\vec{B}_{\K}$.

\subsection*{The quaternionic K\"ahler metric and connection}

Our goal is to match Pedersen and Poon's description of the local geometry of the hyperk\"ahler space ${\cal U}({\cal M})$, centered on the toric symmetries, with Swann's approach, centered on the $\mathbb{H}^{\times}$ structure. The backbone of this correspondence is formed by the equations \eqref{r=qchibq}, \eqref{Higgs-Coulomb} and \eqref{ABC} exhibiting the $q$-dependence of the monopole position vectors $\vec{r}^{\,\I}$, Higgs fields $\Phi_{\I\J}$ and connection 1-forms $A_{\I}$. 

To connect the generalized Gibbons-Hawking form \eqref{GH-Ansatz} with the form \eqref{G-sigma-basis} of the metric we need to evaluate additionally the scalar product $d\vec{r}^{\,\I} \!\cdot d\vec{r}^{\,\J}$ in the inhomogeneous sigma basis. From \eqref{r=qchibq} we get:
\begin{align}
d\vec{r}^{\,\I} \!\sprod d\vec{r}^{\,\J} = |q|^4[&\, d\vec{\chi}^{\,\I} \!\!\cdot\! d\vec{\chi}^{\,\J} + 2\, d(\vec{\chi}^{\,\I} \!\!\cdot\! \vec{\chi}^{\,\J})\, \sigma_0 + 4\, \vec{\chi}^{\,\I} \!\!\cdot\! \vec{\chi}^{\,\J}\sigma_0^2  \nonumber \\[4pt] 
- & 2\,\epsilon_{ijk\,}(\chi^{\I}_i d \chi^{\J}_j - \chi^{\J}_j d \chi^{\I}_i)\,\sigma_k + 4(\vec{\chi}^{\,\I} \!\!\cdot\! \vec{\chi}^{\,\J}\delta_{ij} - \chi^{\I}_i \chi^{\J}_j)\, \sigma_i\sigma_j] \label{dr.dr}
\end{align}
Clearly, the only terms in $G$ that, upon substitution, will contain $\sigma_0$ come from the first term in \eqref{GH-Ansatz}. Specifically, these are
\begin{equation}
|q|^2[U_{\I\J}d(\vec{\chi}^{\,\I} \!\!\cdot\! \vec{\chi}^{\,\J})\pt\sigma_0 + 2 U_{\I\J}\vec{\chi}^{\,\I} \!\!\cdot\! \vec{\chi}^{\,\J} \sigma_0^2] = |q|^2 [d\V\sigma_0 + \V\sigma_0^2]
\end{equation}
and it is a simple exercise to show that they can be cast in the form on the right hand side. Notice though that the term linear in $\sigma_0$ has no correspondent in \eqref{G-sigma-basis}! This problem can be circumvented by an appropriate rescaling of $q$. Under a redefinition $q \to \lambda q$ with $\lambda > 0$, the $\sigma_i$ remain invariant while $\sigma_0 \to \sigma_0 + d \ln \lambda$. By choosing in particular $\lambda = 1/\sqrt{V}$,  the term linear in $\sigma_0$ is eliminated:
\begin{equation}
|q|^2 [d\V\sigma_0 + \V\sigma_0^2] \to |q|^2 \bigg[\sigma_0^2 - \left(\frac{d\V}{2\V}\right)^{\!\!2} \bigg]
\end{equation}

Of course, all the results that we have obtained so far will have to be adjusted to this new scale. In view of the above observation with respect to the scaling behaviour of the left $\mathbb{H}^{\times}$-invariant 1-forms, this is done very easily.  Thus, the monopole position equation \eqref{r=qchibq} becomes
\begin{equation}
\vec{r}^{\,\I} = \frac{q\,\vec{\chi}^{\,\I}\bar{q}}\V \label{scale-pos}
\end{equation}
the Higgs field equation \eqref{Higgs-Coulomb} changes to
\begin{equation}
\Phi_{\I\J} = \frac{\V U_{\I\J}}{|q|^2}
\end{equation}
while, on the other hand, the connection 1-form equation \eqref{ABC} remains unaltered.

We are now fully equipped to bridge the final link between the two descriptions of ${\cal U}({\cal M})$. By substituting these relations as well as the rescaled version of equation \eqref{dr.dr} into the generalized Gibbons-Hawking Ansatz \eqref{GH-Ansatz} we arrive, after some simple formal manipulations, to the hyperk\"ahler cone form \eqref{G-HKC}. In the process we obtain explicit expressions for the $Sp(1)$ connection $\vec{\omega}$ and the metric $g$ of the quaternionic K\"ahler manifold ${\cal M}$. These read as follows:
\begin{equation}
\boxed{ \vec{\omega} = \frac{1}\V [(d\psi_{\I}+C_{\I})\vec{\chi}^{\,\I} - U_{\I\J}\,\vec{\chi}^{\,\I} \!\times d\vec{\chi}^{\,\J}] } \label{Sp(1)-conn}
\end{equation}
and 
\begin{equation}
\boxed{ s\pt g = \frac{1}{2\V}[U_{\I\J}\,d\vec{\chi}^{\,\I} \!\!\cdot\! d\vec{\chi}^{\,\J} + U^{\I\J}(d\psi_{\I}+C_{\I})(d\psi_{\J}+C_{\J})] - |\boldsymbol{\omega}|^2 } \label{qK-metric}
\end{equation}
respectively. By definition, $\boldsymbol{\omega} = \omega_0 + \omega$ is a quaternionic-valued $1$-form with $\mbox{Im\,}\mathbb{H}$-valued part $\omega$ given, in vector form, by \eqref{Sp(1)-conn} and $\mbox{Re\,}\mathbb{H}$-valued part $\omega_0$ given by
\begin{equation}
\omega_0 = \frac{d\V}{2\V} \label{omega_0}
\end{equation}
This $4n$-dimensional metric has $n+1$ manifest commuting Killing vector fields given by $X^{\I}=\partial_{\psi_I}$. As one can easily check, the $X^{\I}$ preserve the $Sp(1)$ connection as well, which means that the the corresponding vector-valued functions $\vec{R}$ that appear (in quaternionic form) in the  equations \eqref{R-eq} and \eqref{QK-rot} are all vanishing in this case.

\section{The case of four dimensions} \label{SEC:4d}

In this section we will focus our attention to the $n=1$ case, which is the case of self-dual Einstein manifolds with two linearly independent commuting Killing vector fields. The local geometry of these manifolds has been characterized and completely classified by Calderbank and Pedersen in the remarkable paper \cite{MR1950174} through an essentially intrinsic, four-dimensional approach. We will show that our results, when particularized to $n=1$, reproduce those of Calderbank and Pedersen and thus add a new perspective to the problem.

For $n=1$ the indices $I,J,K$ take only two values, 0 and 1. There are two monopoles, with position vectors
\begin{align}
\vec{r}^{\,\0} & = \frac{q\pt \k\pt \bar{q} \label{r0}}{\V} \\[4pt]
\vec{r}^{\,\1} & = \frac{q\pt(\rho\pt\j + \eta\pt\k)\pt\bar{q}}{\V}
\end{align}
For simplicity, we drop the index $1$ from $\rho^{\1}$ and $\eta^{\1}$, since these are the only coordinates of this kind and there is no danger of confusion.

The $n=1$ case is degenerate, in the sense that many of the equations that we wrote for generic $n$ become trivial or do not apply for $n=1$, so we are left in effect with a smaller number of relations. Thus, from the set of equations \eqref{curlxB1}--\eqref{curlxB3} only equation \eqref{curlxB2} survives and reads
\begin{equation}
 \partial_{\rho}{(B_{\K})}_3 = \partial_{\eta}{(B_{\K})}_2 
\end{equation}
The equations \eqref{dchiB}--\eqref{detaB} have either trivial or no content. The equation \eqref{partial-tr} can be re-written in this case as follows
\begin{equation}
\partial_{\rho}{(B_{\K})}_2 + \partial_{\eta}{(B_{\K})}_3 = {(B_{\K})}_2/\rho
\end{equation}
These last two relations arise in the context of \cite{MR1950174} as consequences of a so-called Joyce equation. Furthermore, equations \eqref{B1_1}--\eqref{B0_3} and \eqref{B1=d1V}--\eqref{B3=d3V} give, respectively,  
\begin{align}
{(B_{\0})}_1 & = 0 & {(B_{\1})}_1 & = 0 \nonumber \\[4pt]
{(B_{\0})}_2 & = \rho \V_{\eta} - \eta \V_{\rho}  & {(B_{\1})}_2 & = \V_{\rho} \nonumber \\[4pt]
{(B_{\0})}_3 & = \V - \rho \V_{\rho} - \eta \V_{\eta}   & {(B_{\1})}_3 & = \V_{\eta}
\end{align}
The equations \eqref{ddV1}--\eqref{ddV3} yield no constraints at all on $\V$. On the other hand, the whole set of second-order differential constraints \eqref{dubya1}--\eqref{dubya3} reduces to a single one, namely
\begin{equation}
\rho(\V_{\rho\rho} + \V_{\eta\eta}) = \V_{\rho}
\end{equation}
corresponding to equation \eqref{dubya2} with $I=1$ and $j=3$. By substituting $\V = \sqrt{\rho}\, \V'$, this constraint can be reformulated as an eigenvalue equation for the Laplacian on the upper half-plane endowed with the Poincar\'e metric, which features prominently in Calderbank and Pedersen's paper.

The reduced Bogomol'nyi equation \eqref{Bogo} reads simply 
\begin{equation}
dC_{\K} = 0 \label{Bogo-4d}
\end{equation}
Locally at least, we can always pick $C_{\K} = 0$. The Higgs field components $U_{\I\J}$ take the simplest form when expressed in the modified frame 
\begin{equation}
\left(\! \begin{array}{c}
\alpha \\
\beta
\end{array} \!\right)
=
\left( \begin{array}{cc} 
\rho & 0 \\ 
\eta & 1
\end{array} \right) 
\!
\left(\! \begin{array}{c}
d\psi_{\1} \\
d\psi_{\0}
\end{array} \!\right)
\end{equation}
that we have introduced in connection with equations \eqref{UV}. From the latter we have
\begin{equation}
\left( \begin{array}{cc} 
\rho & 0 \\ 
\eta & 1
\end{array} \right) 
\!
\left(\! \begin{array}{cc} 
U_{\1\1} & U_{\1\0} \\ 
U_{\0\1} & U_{\0\0} 
\end{array} \!\right)
\!
\left( \begin{array}{cc} 
\rho & \eta  \\ 
0      &  1
\end{array} \right) 
= \frac{1}{2}
\left( \begin{array}{cc} 
\rho \V_{\rho} & \rho \V_{\eta} \\ 
\rho \V_{\eta} & \V - \rho \V_{\rho} 
\end{array} \!\right) 
\end{equation}
By substituting these expressions into the equations \eqref{Sp(1)-conn}, \eqref{qK-metric} and \eqref{omega_0} we obtain directly
\begin{equation}
g = \frac{\V\V_{\rho}-\rho \V_{\rho}^2-\rho \V_{\eta}^2}{4\rho \V^2} (d\rho^2+d\eta^2) + \frac{[(\V-\rho \V_{\rho})\alpha - \rho \V_{\eta}\beta]^2 + [\rho \V_{\eta}\alpha - \rho \V_{\rho}\beta]^2}{\rho \V^2(\V\V_{\rho}-\rho \V_{\rho}^2-\rho \V_{\eta}^2)}
\end{equation}
and
\begin{equation}
\boldsymbol{\omega} = \frac{\V_{\rho}d\rho + \V_{\eta}d\eta}{2\V} + \frac{\V_{\eta}d\rho - \V_{\rho}d\eta}{2\V} \i + \frac{\alpha}{\V} \j + \frac{\beta}{\V} \k
\end{equation}
As one can easily verify, these are indeed, up to simple redefinitions, the Calderbank-Pedersen formulas for the quaternionic K\"ahler metric and $Sp(1)$ connection 1-form (strictly speaking, only the purely imaginary part of $\boldsymbol{\omega}$ corresponds to the connection 1-form).

\section{Explicit quaternionic K\"ahler metrics from \\ the generalized Legendre transform}

We will now show how to construct explicit solutions of the reduced Bogomol'nyi equation \eqref{Bogo} as well as of the second-order differential constraints \eqref{ddV1}--\eqref{ddV3} and \eqref{dubya1}--\eqref{dubya3} on $\V$ by means of the generalized Legendre transform procedure, for a large class of quaternionic K\"ahler manifolds, namely the ones sitting in the image of the so-called c-map. 

The c-map, or more precisely, the local c-map, associates to any projective special K\"ahler manifold of complex dimension $n-1$ a quaternionic K\"ahler manifold of quaternionic dimension $n$ with $2n+6$ isometries, of which $n+1$ are commuting \cite{Cecotti:1988qn,Ferrara:1989ik}. In physics, where it was originally discovered, it corresponds to the dimensional reduction of $N=2$ supergravity from four to three space-time dimensions, followed by the dualization of the vector multiplets into hypermultiplets. Consider two such manifolds related by the c-map. Over the projective special K\"ahler manifold one can construct a canonical $\mathbb{C}^{\times}$-bundle whose total space has a special K\"ahler structure in the affine sense (we use here the terminology of \cite{Freed:1997dp}). The local geometry of the projective special K\"ahler base is completely determined by a holomorphic function ${\cal F}$ on this space -- the prepotential, homogeneous of degree two in its $n$ variables. On the other hand, over the quaternion K\"ahler manifold one can build a canonical $\mathbb{H}^{\times}$-bundle -- the Swann bundle, whose total space carries a hyperk\"ahler structure. The local geometry of the quaternion K\"ahler base is encoded in the meromorphic function ${\cal H}$ of the generalized Legendre transform construction, homogeneous of degree one in its $n+1$ variables. In \cite{Rocek:2005ij,Rocek:2006xb} Ro\v{c}ek, Vafa and Vandoren found that these two functions are related to one another in a remarkably simple and direct way, namely
\begin{equation}
{\cal H}(\hat{\eta}^{\0},\hat{\eta}^{\1},\cdots,\hat{\eta}^{n}) = \frac{{\cal F}(\hat{\eta}^1,\cdots,\hat{\eta}^{n})}{\hat{\eta}^0} \label{c-map}
\end{equation}
The contour $\Gamma$ of \eqref{GLT-action} is in this case an eight-shape curve around the two roots of $\hat{\eta}^0$. 

Given a prepotential ${\cal F}$, the contour integral $F$ can be calculated explicitly by means of the residue theorem in terms of the parameters of the ${\cal O}(2)$ sections $\hat{\eta}^{\I}$. One has 
\begin{equation}
F = \frac{{\cal F}(\hat{\eta}^{\1}(\zeta^0_+),\cdots,\hat{\eta}^{n}(\zeta^0_+)) + {\cal F}(\hat{\eta}^{\1}(\zeta^0_-),\cdots,\hat{\eta}^{n}(\zeta^0_-))}{r^{\0}}
\end{equation}
where $\zeta^0_{\pm} = (x^{\0} \pm r^{\0})/ 2z^{\0}$ are the antipodally-conjugated roots of $\hat{\eta}^0$. The ensuing steps are standard generalized Legendre transform procedure: one computes further its first and second derivatives and then the Higgs field components $\Phi_{\I\J}$, connection 1-forms $A_{\K}$ and hyperk\"ahler potential ${\cal K}$.

To descend onto the quaternionic K\"ahler base we must substitute the pivotal relations \eqref{r_0}--\eqref{r_I} into the results. In practice, this works well for $\Phi_{\I\J}$ and ${\cal K}$: it yields a $\Phi_{\I\J}$ of the form \eqref{Higgs-Coulomb}, therefore allowing us to identify $U_{\I\J}$, and a ${\cal K}$ of the form
\begin{equation}
{\cal K} = \V |q|^2 \label{hyperK-pot}
\end{equation}
therefore allowing us to quickly identify $\V$. Indeed, based on collective rotation and scaling behavior arguments we know that ${\cal K}$ has to be proportional to $|q|^2$; the proportionality factor must be $\V$, so that at the rescaling of $q$ mandated in the last subsection of section \ref{SEC:H*-red}, ${\cal K}$ becomes eventually equal to $|q|^2$. On the other hand, when it comes to the connection 1-forms $A_{\K}$, one runs into problems. By analyzing concrete examples we found that the result of substituting \eqref{r_0}--\eqref{r_I} into \eqref{A-GLT} cannot be cast into the form \eqref{ABC}. Substituting instead the scaled version \eqref{scale-pos} does not make a difference. Neither does choosing different reference configurations \eqref{r_0}--\eqref{r_I}. This discrepancy should give us cause for serious concern, for it signals that if we substitute the outcome further, along with the previously obtained expressions for $\Phi_{\I\J}$ and ${\cal K}$, into the generalized Gibbons-Hawking Ansatz \eqref{GH-Ansatz}, we will fail to come across the Swann form of the metric.

The key to understanding and eventually circumventing these difficulties lays in the realization that the relevant quantity is not $A_{\K}$ {\it per se} but the cohomology class of $A_{\K}$ and that thus we have at our disposal a certain freedom in choosing the particular representative of the equivalence class. The form \eqref{A-GLT} may simply not be a convenient representative  as far as taking the Swann quotient is concerned. We propose instead a different, canonical choice of representative that is not riddled with these deficiencies and delivers the desired outcome. It is obtained by shifting the representative form \eqref{A-GLT} as follows:
\begin{equation}
\boxed{ A_{\K} \to A_{\K} - d\psi^{\prime}_{\K} \qquad\mbox{with}\qquad \psi^{\prime}_{\K} = c_{\K}\mbox{\pt Im\pt}(z^{\0}F_{z^0x^K}) } \label{shifted-A}
\end{equation}
where $c_{\K} = 1$ for $K \neq 0$ and $c_{\0}=1/2$. The evidence supporting this claim comes form a host of particular examples which we will discuss  in detail shortly.

What we have described here is in essence a quotienting procedure: one starts with a Swann bundle characterized by a meromorphic function of the form \eqref{c-map} and ends up with the local geometry ({\it i.e.} metric, $Sp(1)$ connection) of a quaternionic K\"ahler manifold.  At the core of this approach there are two gauge choices: that of a reference configuration -- the equations \eqref{r_0}--\eqref{r_I}, and that of a cohomology class representative -- the equation \eqref{shifted-A}.

\subsection*{A $2\times{\cal O}(2)$ model}

We start by reviewing the simplest example of such a construction, which has been discussed from a somehow different perspective in \cite{Anguelova:2004sj} in relation to the universal hypermultiplet of string theory. The generating $F$-function is based on the holomorphic prepotential ${\cal F}(X_1) = X_1^2$ and reads\footnote{Here and in the remainder of these notes we alter our previous notation conventions and lower all indices in order to avoid an excessive use of parentheses that would otherwise be required to prevent  possible confusions with exponents.}
\begin{equation}
F = \frac{1}{2\pi i} \oint_{\Gamma_0} \frac{d\zeta}{\zeta}\, \frac{(\hat{\eta}_1)^2}{\hat{\eta}_0}
\end{equation}
According to \cite{Cecotti:1988qn}, it gives the pseudo-quaternionic K\"ahler metric on the symmetric space $SU(2,1)/S(U(2) \times U(1))$.

A standard generalized Legendre transform calculation produces the following hyperk\"ahler potential 
\begin{equation}
{\cal K} = \frac{2(\vec{r}_0\vprod\vec{r}_1)^2}{r_0^3}
\end{equation}
and Higgs field components 
\begin{align}
\Phi_{\0\0} & = \frac{3(\vec{r}_0\sprod\vec{r}_1)^2-r_0^2r_1^2}{r_0^5} \nonumber \\
\Phi_{\0\1} & = - \frac{2\,\vec{r}_0\sprod\vec{r}_1}{r_0^3} \nonumber \\
\Phi_{\1\1} & = \frac{2}{r_0}
\end{align} 
The way we wrote them, these expressions are manifestly invariant at collective rotations of the monopole configuration and are homogeneous of degree 1 respectively -1 at a simultaneous scaling of the monopole position vectors. Following the prescription stated above and substituting for $\vec{r}_0$ and $\vec{r}_1$ the formulas \eqref{r_0} and \eqref{r_1}, we obtain that ${\cal K}$ is of the form \eqref{hyperK-pot} with
\begin{equation}
\V = 2\rho^2
\end{equation}
while $\Phi_{\I\J}$ is of the form \eqref{Higgs-Coulomb} with
\begin{align}
U_{\0\0} & = 2\eta^2-\rho^2 \nonumber \\[4pt]
U_{\0\1} & = -2\eta \nonumber \\[4pt]
U_{\1\1} & = 2
\end{align}
Just like in section \ref{SEC:4d}, we dropped for simplicity the index $1$ from $\rho^{\1}$ and $\eta^{\1}$. Using the definition \eqref{B-def} we get further
\begin{align}
{(B_{\0})}_1 & = 0 & {(B_{\1})}_1 & = 0 \nonumber \\[4pt]
{(B_{\0})}_2 & = -4\eta\rho & {(B_{\1})}_2 & = 4\rho \nonumber \\[4pt]
{(B_{\0})}_3 & = -2\rho^2 & {(B_{\1})}_3 & = 0 \label{Bees}
\end{align} 
The connection 1-forms $A_{\K}$ can be also calculated explicitly in terms of the components of the vectors $\vec{r}_0$ and $\vec{r}_1$. One can check directly that if one chooses for the cohomology class of $A_{\K}$ the representative \eqref{A-GLT} one does not obtain a result of the form \eqref{ABC} upon substitution of the expressions \eqref{r_0} and \eqref{r_1} for $\vec{r}_0$ and $\vec{r}_1$. If, on the other hand, one picks the shifted representative \eqref{shifted-A}, one does get the form \eqref{ABC}, with the $\vec{B}_{\K}$ coefficients given precisely by \eqref{Bees} and
\begin{equation}
C_{\0} = C_{\1} = 0
\end{equation}
This is trivially a solution to the reduced Bogomol'nyi equation \eqref{Bogo-4d}.

\subsection*{A $3\times{\cal O}(2)$ model}

We apply next the same procedure to the Swann bundle with generating function 
\begin{equation}
F = \frac{1}{2\pi i} \oint_{\Gamma_0} \frac{d\zeta}{\zeta}\, \frac{\hat{\eta}_1\hat{\eta}_2}{\hat{\eta}_0}
\end{equation} 
based on the holomorphic prepotential ${\cal F}(X_1,X_2)=X_1X_2$. The quotient construction returns this time the eight-dimensional symmetric pseudo-quaternionic K\"ahler geometry of $SU(2,2)/S(U(2) \times U(2))$.

The generalized Legendre transform of $F$ yields the hyperk\"ahler potential
\begin{equation}
{\cal K} = \frac{2(\vec{r}_0\vprod\vec{r}_1)\sprod(\vec{r}_0\vprod\vec{r}_2)}{r_0^3} 
\end{equation}
and the Higgs field components
\begin{align}
\Phi_{\0\0} & = \frac{3(\vec{r}_0\sprod\vec{r}_1)(\vec{r}_0\sprod\vec{r}_2)-r_0^2(\vec{r}_1\sprod\vec{r}_2)}{r_0^5} \hspace{-50pt}& \Phi_{\0\1} & = -\frac{\vec{r}_0\sprod\vec{r}_2}{r_0^3} \nonumber \\
\Phi_{\1\1} & = 0 & \Phi_{\0\2} & = - \frac{\vec{r}_0\sprod\vec{r}_1}{r_0^3} \nonumber \\
\Phi_{\2\2} & = 0 & \Phi_{\1\2} & = \frac{1}{r_0} 
\end{align}
Substitution of the corresponding formulas \eqref{r_0}--\eqref{r_I} for $\vec{r}_0$, $\vec{r}_1$ and $\vec{r}_2$ gives
\begin{equation}
\V = 2\rho_1\rho_2 \label{Vee}
\end{equation}
respectively
\begin{align}
U_{\0\0} & = 2\eta_1\eta_2 - \rho_1\rho_2 \hspace{-50pt} & U_{\0\1} & = -\eta_2 \nonumber \\[4pt]
U_{\1\1} & = 0 & U_{\0\2} & = -\eta_1 \nonumber  \\[4pt]
U_{\2\2} & = 0 & U_{\1\2} & = 1
\end{align}
From \eqref{B-def} one then gets the following $\vec{B}_{\K}$ field components:
\begin{align}
{(B_{\0})}_1 & = -2\eta_1\chi_2 & {(B_{\1})}_1 & = 2\chi_2  & {(B_{\2})}_1 & = 0 \nonumber \\[4pt]
{(B_{\0})}_2 & = -2(\eta_1\rho_2+\rho_1\eta_2) & {(B_{\1})}_2 & = 2\rho_2 & {(B_{\2})}_2 & = 2\rho_1 \nonumber \\[4pt]
{(B_{\0})}_3 & = -2\rho_1\rho_2 & {(B_{\1})}_3 & = 0  & {(B_{\2})}_3 & = 0
\end{align}
Substituting the formulas for $\vec{r}_0$, $\vec{r}_1$, $\vec{r}_2$ into the shifted representative 1-forms prescribed by \eqref{shifted-A} returns results of the form \eqref{ABC} with these same $\vec{B}_{\K}$ coefficients and $C_{\K}$ components
\begin{align}
C_{\0} & = \frac{1}{2}(\chi_2d\rho_1-\rho_1d\chi_2) \\[2pt]
C_{\1} & = C_{\2} = 0
\end{align}
One can verify explicitly that these expressions provide indeed a solution to the generalized Bogomol'nyi equation \eqref{Bogo} and that the form \eqref{Vee} of $\V$ satisfies the (applicable) second order differential equations \eqref{ddV1}--\eqref{ddV3} and \eqref{dubya1}--\eqref{dubya3} as well.

\subsection*{Another $3\times{\cal O}(2)$ model}

The previous model is symmetric at the interchange of $\hat{\eta}_1$ and $\hat{\eta}_2$. This accidental symmetry might have potential hidden consequences that could make the gauge choice \eqref{shifted-A} work in this case but not in others. To rule this out, we will consider now a $3 \times {\cal O}(2)$ model which does not have this symmetry. Let
\begin{equation}
F = \frac{1}{2\pi i} \oint_{\Gamma_0} \frac{d\zeta}{\zeta}\, \frac{(\hat{\eta}_2)^3}{\hat{\eta}_0\hat{\eta}_1}
\end{equation}
corresponding to the meromorphic prepotential ${\cal F}(X_1,X_2) = X_2^3/X_1$. According to \cite{Cecotti:1988qn} (for a detailed discussion see also \cite{Gunaydin:2007qq}) this yields the eight-dimensional pseudo-quaternionic K\"ahler geometry of the symmetric space $G_{2(2)}/SO(4)$.

The generalized Legendre transform machinery gives in this case the hyperk\"ahler potential
\begin{equation}
{\cal K} = \frac{2(\vec{r}_0\vprod\vec{r}_1)\sprod(\vec{r}_0\vprod\vec{r}_2)}{r_0^3(\vec{r}_0\vprod\vec{r}_1)^4} [3r_0^2|\vec{r}_0 \sprod (\vec{r}_1\vprod\vec{r}_2)|^2 + |(\vec{r}_0\vprod\vec{r}_1)\sprod(\vec{r}_0\vprod\vec{r}_2)|^2]
\end{equation}
and the Higgs field components
{\allowdisplaybreaks
\begin{align} 
&{\allowdisplaybreaks \begin{aligned}
\Phi_{\0\0} = \frac{1}{r_0^5(\vec{r}_0\vprod\vec{r}_1)^6} \Big( 3r_0^2|\vec{r}_0 \sprod (\vec{r}_1\vprod\vec{r}_2)|^2 [\phantom{3}
& (\vec{r}_0\vprod\vec{r}_1)^2 (\vec{r}_0\sprod\vec{r}_1)(\vec{r}_0\sprod\vec{r}_2) \\[-7pt]
- &  (\vec{r}_0\vprod\vec{r}_1)\sprod(\vec{r}_1\vprod\vec{r}_2) (\vec{r}_0\sprod\vec{r}_1) r_0^2 \\[2pt]
- & (\vec{r}_0\vprod\vec{r}_1)\sprod(\vec{r}_0\vprod\vec{r}_2)  r_0^2r_1^2] \\[2pt]
+ (\vec{r}_0\vprod\vec{r}_1)\sprod(\vec{r}_0\vprod\vec{r}_2) [
 3& (\vec{r}_0\vprod\vec{r}_1)^4(\vec{r}_0\sprod\vec{r}_2)^2  \\[2pt]
- & (\vec{r}_0\vprod\vec{r}_1)^2(\vec{r}_0\vprod\vec{r}_1)\sprod(\vec{r}_1\vprod\vec{r}_2) (\vec{r}_0\sprod\vec{r}_2)r_0^2 \\[2pt]
- & (\vec{r}_0\vprod\vec{r}_1)^2 (\vec{r}_0\vprod\vec{r}_1)\sprod(\vec{r}_0\vprod\vec{r}_2) (\vec{r}_1\sprod\vec{r}_2)r_0^2 \\[-3pt]
+ & 2|(\vec{r}_0\vprod\vec{r}_1)\sprod(\vec{r}_1\vprod\vec{r}_2)|^2r_0^4 ] \Big)
\end{aligned} } \nonumber \\
&\begin{aligned}
\Phi_{\1\1} = - \frac{2(\vec{r}_0\vprod\vec{r}_1)\sprod(\vec{r}_0\vprod\vec{r}_2) }{r_0(\vec{r}_0\vprod\vec{r}_1)^6} \Big( 3r_0^2|\vec{r}_0 \sprod (\vec{r}_1\vprod\vec{r}_2)|^2-|(\vec{r}_0\vprod\vec{r}_1)\sprod(\vec{r}_0\vprod\vec{r}_2)|^2 \Big)
\end{aligned} \nonumber \\
&\begin{aligned}
\Phi_{\2\2} =  \frac{6(\vec{r}_0\vprod\vec{r}_1)\sprod(\vec{r}_0\vprod\vec{r}_2) }{r_0(\vec{r}_0\vprod\vec{r}_1)^2}
\end{aligned} \nonumber  \\
&\begin{aligned}
\Phi_{\0\1} = \frac{1}{r_0^3(\vec{r}_0\vprod\vec{r}_1)^6} 
\Big( 3r_0^2|\vec{r}_0 \sprod (\vec{r}_1\vprod\vec{r}_2)|^2 [
   2& (\vec{r}_0\vprod\vec{r}_1)\sprod(\vec{r}_1\vprod\vec{r}_2)r_0^2 
+ (\vec{r}_0\vprod\vec{r}_1)^2(\vec{r}_0\sprod\vec{r}_2)]  \\[-6pt]
- |(\vec{r}_0\vprod\vec{r}_1)\sprod(\vec{r}_0\vprod\vec{r}_2)|^2 [
  2 & (\vec{r}_0\vprod\vec{r}_1)\sprod(\vec{r}_1\vprod\vec{r}_2)r_0^2 
-  (\vec{r}_0\vprod\vec{r}_1)^2(\vec{r}_0\sprod\vec{r}_2)] \Big)
\end{aligned} \nonumber \\
&\begin{aligned}
\Phi_{\0\2} = -\frac{3}{r_0^3(\vec{r}_0\vprod\vec{r}_1)^4} \Big( r_0^2|\vec{r}_0 \sprod (\vec{r}_1\vprod\vec{r}_2)|^2  (&\vec{r}_0\sprod\vec{r}_1)   \nonumber \\[-7pt]
\! - (\vec{r}_0\vprod\vec{r}_1)\sprod(\vec{r}_0\vprod\vec{r}_2)[&
   (\vec{r}_0\vprod\vec{r}_1)\sprod(\vec{r}_1\vprod\vec{r}_2)r_0^2  
-  (\vec{r}_0\vprod\vec{r}_1)^2(\vec{r}_0\sprod\vec{r}_2)] \Big) 
\end{aligned} \nonumber \\
&\begin{aligned}
\Phi_{\1\2} = \frac{3}{r_0(\vec{r}_0\vprod\vec{r}_1)^4} \Big( r_0^2|\vec{r}_0 \sprod (\vec{r}_1\vprod\vec{r}_2)|^2 - |(\vec{r}_0\vprod\vec{r}_1)\sprod(\vec{r}_0\vprod\vec{r}_2)|^2 \Big)
\end{aligned} 
\end{align} 
}These are rather complicated but on the other hand manifestly rotation-invariant and homogeneous expressions. It takes a certain amount of effort to arrange them in this way, but keep in mind that in practice one need not do so. We did it here mainly to show that it is possible and to maintain at the same time a parallel narrative with the previous examples. Upon substituting the expressions \eqref{r_0}--\eqref{r_I} for $\vec{r}_0$, $\vec{r}_1$, $\vec{r}_2$, the formulas become sensibly simpler. One gets:
\begin{equation}
\V = \frac{2\rho_2(\rho_2^2+3\chi_2^2)}{\rho_1}
\end{equation}
respectively
{\allowdisplaybreaks
\begin{align}
U_{\0\0} & = \frac{2\eta_1^2\rho_2(\rho_2^2-3\chi_2^2)}{\rho_1^3} - \frac{6\eta_1\eta_2(\rho_2^2-\chi_2^2)}{\rho_1^2} + \frac{\rho_2(6\eta_2^2-\rho_2^2-3\chi_2^2)}{\rho_1} \nonumber \\
U_{\1\1} & = \frac{2\rho_2(\rho_2^2-3\chi_2^2)}{\rho_1^3} \nonumber \\
U_{\2\2} & = \frac{6\rho_2}{\rho_1} \nonumber \\
U_{\0\1} & = - \frac{2\eta_1\rho_2(\rho_2^2-3\chi_2^2)}{\rho_1^3} + \frac{3\eta_2(\rho_2^2-\chi_2^2)}{\rho_1^2} \nonumber \\
U_{\0\2} & = \frac{3\eta_1(\rho_2^2-\chi_2^2)}{\rho_1^2} - \frac{6\eta_2\rho_2}{\rho_1} \nonumber \\
U_{\1\2} & = - \frac{3(\rho_2^2-\chi_2^2)}{\rho_1^2} 
\end{align} 
}As shown earlier, the fields $U_{\I\J}$ can in fact be derived from the function $\V$ alone, which encodes all local geometric information. One can check directly that this $\V$ satisfies the differential constraints \eqref{ddV1}--\eqref{ddV3} and \eqref{dubya1}--\eqref{dubya3}. The $\vec{B}_{\K}$ field components read
{\allowdisplaybreaks
\begin{align*}
{(B_{\0})}_1\! & =\! \frac{6\eta_1\chi_2(\rho_2^2-\chi_2^2)}{\rho_1^2} \!-\! \frac{12\eta_2\rho_2\chi_2}{\rho_1} & \!\!\!{(B_{\1})}_1\! & =\! - \frac{6\chi_2(\rho_2^2-\chi_2^2)}{\rho_1^2} & \!\!\!{(B_{\2})}_1\! & =\!  \frac{12\rho_2\chi_2}{\rho_1} \\
{(B_{\0})}_2\! & =\! \frac{2\eta_1\rho_2(\rho_2^2+3\chi_2^2)}{\rho_1^2} \!-\! \frac{6\eta_2(\rho_2^2+\chi_2^2)}{\rho_1} & \!\!\!{(B_{\1})}_2\! & =\! - \frac{2\rho_2(\rho_2^2+3\chi_2^2)}{\rho_1^2}  & \!\!\!{(B_{\2})}_2\! & =\! \frac{6(\rho_2^2+\chi_2^2)}{\rho_1} \\
{(B_{\0})}_3\! & =\! - \frac{2\rho_2(\rho_2^2+3\chi_2^2)}{\rho_1}  & \!\!\!{(B_{\1})}_3\! & = 0 & \!\!\!{(B_{\2})}_3\! & = 0 
\end{align*} 
}Just as before, substituting the formulas \eqref{r_0}--\eqref{r_I} for $\vec{r}_0$, $\vec{r}_1$, $\vec{r}_2$ into the modified representative \eqref{shifted-A} leads to generalized Gibbons-Hawking connection 1-forms of the form \eqref{ABC} with $\vec{B}_{\K}$ given above and the following sigma-free components:
{\allowdisplaybreaks
\begin{align}
C_{\0} {} = & \phantom{+\ } \frac{3\chi_2}{2\rho_1^4}(\rho_1^2\chi_2^2-\rho_1^2\rho_2^2-2\rho_1^2\eta_2^2+2\eta_1^2\chi_2^2-6\eta_1^2\rho_2^2+8\rho_1\eta_1\rho_2\eta_2)d\rho_1 \nonumber \\
& - \frac{3}{2\rho_1^3} (\rho_1^2\chi_2^2+\rho_1^2\rho_2^2-2\rho_1^2\eta_2^2+2\eta_1^2\chi_2^2-2\eta_1^2\rho_2^2+4\rho_1\eta_1\rho_2\eta_2)d\chi_2 \nonumber \\
& + \frac{3\chi_2}{\rho_1^3}(\rho_1^2\rho_2+2\eta_1^2\rho_2-2\rho_1\eta_1\eta_2)d\rho_2 \nonumber \\
C_{\1} {} = & - \frac{6\chi_2}{\rho_1^4}(\eta_1\chi_2^2-3\eta_1\rho_2^2+2\rho_1\rho_2\eta_2)d\rho_1 
+  \frac{6}{\rho_1^3}\ (\eta_1\chi_2^2 -\eta_1\rho_2^2 + \rho_1\rho_2\eta_2)d\chi_2 \nonumber \\
& -  \frac{6\chi_2}{\rho_1^3}(2\eta_1\rho_2-\rho_1\eta_2)d\rho_2 \nonumber \\
C_{\2} {} = & - \frac{6\chi_2}{\rho_1^3}(2\eta_1\rho_2-\rho_1\eta_2)d\rho_1 + \frac{6}{\rho_1^2}(\eta_1\rho_2-\rho_1\eta_2)d\chi_2 + \frac{6\chi_2}{\rho_1^2}\eta_1d\rho_2
\end{align} 
}These expressions constitute yet another explicit and non-trivial solution to the reduced Bogomol'nyi equation \eqref{Bogo}.

\subsection*{A twisted look at flat space}

We end up with a discussion of a toy model based on a single ${\cal O}(2)$ section which does not fit into the c-map-related class of examples that we have focused on so far, but which nevertheless shares some interesting features with these. Consider the function \cite{Karlhede:1984vr}
\begin{equation}
F = \frac{1}{2\pi i} \oint_{\Gamma_0} \frac{d\zeta}{\zeta}\, \hat{\eta}_0 \ln \hat{\eta}_0
\end{equation}
with the contour $\Gamma_0$ surrounding the two roots of $\hat{\eta}_0$ over the logarithmic branch cuts. The contour integral can be evaluated explicitly in terms of the parameters of the section $\hat{\eta}_0$:
\begin{equation}
F = r_0 - x_0\, \mbox{arctanh}\frac{x_0}{r_0}
\end{equation}
We use the notations introduced in \eqref{O2-multiplet} but with the indices lowered. The Legendre transform of $F$ yields
\begin{equation}
{\cal K} = r_0
\end{equation}
The equations \eqref{Phi-GLT} and \eqref{A-GLT} give in turn
\begin{align}
\Phi_{\0\0} & = \frac{1}{2r_0} \\
A_{\0}\ & = -\frac{x_0}{2r_0} \mbox{\,Im\,}\frac{dz_0}{z_0}
\end{align}
If we now substitute in the last equation the expression \eqref{r_0} for $\vec{r}_0 = 2\mbox{\pt Im\,}z_0\pt\i - 2\mbox{\pt Re\,}z_0\pt\j + x_0\pt\k$ we obtain
\begin{equation}
A_{\0} = \frac{q_0^2+q_3^2-q_1^2-q_2^2}{2(q_0^2+q_3^2)(q_1^2+q_2^2)} [(q_1q_3-q_0q_2)\sigma_1 + (q_2q_3+q_0q_1)\sigma_2]
\end{equation}
Just like in the previous examples, we find that the cohomology class representative \eqref{A-GLT} does not deliver a result of the form \eqref{ABC}. But observe that this can be written as follows
\begin{equation}
A_{\0} = \sigma_3 + d\psi^{\prime}_{\0} \qquad\mbox{with}\qquad \psi^{\prime}_{\0} = \frac{1}{2}\arctan\frac{q_0q_1+q_2q_3}{q_0q_2-q_1q_3}
\end{equation}
On the other hand, a simple counting argument implies that in four dimensions, which is our case, the toric coordinate $\psi_0$ defined by equation \eqref{psi-coord} cannot be independent from the $q$ variables: were it independent, one would then describe a four-dimensional space by means of five free parameters! We stress that this is a singular occurrence: for models based on two or more ${\cal O}(2)$ sections the argument clearly fails and the contrary statement holds.  We must have therefore $\psi_{\0} = \psi/2$, with $\psi$  the Euler variable defined in \eqref{Euler-coords}, and so
\begin{equation}
\psi_{\0} = \frac{1}{2}\arctan\frac{2(q_1q_2-q_0q_3)}{q_0^2+q_1^2-q_2^2-q_3^2}
\end{equation}
These two observations suggest we choose a different representative for the cohomology class of $A_{\0}$ as follows
\begin{equation}
A_{\0} \to A_{\0} - d(\psi_{\0}+\psi^{\prime}_{\0})
\end{equation}
The shift term can be expressed in terms of the original variables. By way of the $\arctan$ addition theorem we have
\begin{equation}
\psi_{\0}+\psi^{\prime}_{\0} = \frac{1}{2} \arctan\frac{r_0\mbox{\,Re\,}z_0}{x_0\mbox{\,Im\,}z_0}
\end{equation}
Using equation \eqref{r_0} one can show that that $r_0 = |q|^2$ and $d\vec{r}_0 \cdot d\vec{r}_0= 4|q|^4(\sigma_0^2+\sigma_1^2+\sigma_2^2)$. Substituting everything into the Gibbons-Hawking formula, we  get finally
\begin{equation}
G = |q|^2(\sigma_0^2+\sigma_1^2+\sigma_2^2+\sigma_3^2) = dr^2 + r^2(\sigma_1^2+\sigma_2^2+\sigma_3^2)
\end{equation}
that is, the flat metric on $\mathbb{R}^4 \!\smallsetminus\! \{{\mathbf 0}\}$.

\vskip30pt

\subsection*{\sf Acknowledgements}

The author wishes to thank Martin Ro\v{c}ek for his constant support and enlightening discussions.  This paper has also benefited from the useful comments of Boris Pioline and Stefan Vandoren. Many thanks go to Maria Victoria Fernandez-Serra and Moustapha Thioye for providing critical computer muscle, without which one of the computations presented here could not have been performed.

\vskip50pt

\bibliographystyle{utcaps}
\bibliography{ToricHKC}

\end{document}